\documentclass[12pt]{elsarticle}

\usepackage{geometry} 
\usepackage{color} 
\usepackage{graphicx}
\usepackage{caption}
\usepackage{subcaption}
\usepackage{courier}
\usepackage{amsmath, amssymb}
\usepackage{setspace}
\usepackage{pdflscape}

\date{} 
\usepackage[toc,page]{appendix}

\begin{document}
\begin{frontmatter}

\title{Time-dependent stochastic basis adaptation for uncertainty quantification}
 
\author{Ramakrishna Tipireddy}
\ead{ramakrishna.tipireddy@pnnl.gov}

\author{Panos Stinis}
\address{Pacific Northwest National Laboratory, Richland, WA}

\author{Alexandre M. Tartakovsky}
\address{Pacific Northwest National Laboratory, Richland, WA; Department of Civil and Environmental Engineering, University of Illinois Urbana-Champaign, Urbana, IL}

\begin{abstract}
We extend stochastic basis adaptation and spatial domain decomposition methods to solve time varying stochastic partial differential equations (SPDEs) with a large number of input random parameters. Stochastic basis adaptation allows the determination of a low dimensional stochastic basis representation of a quantity of interest (QoI). Extending basis adaptation to time-dependent problems is challenging because small errors introduced in the previous time steps of the low dimensional approximate solution accumulate over time and cause divergence from the true solution. To address this issue we have introduced an approach where the basis adaptation varies at every time step so that the low dimensional basis is adapted to the QoI at that time step. We have coupled the time-dependent basis adaptation with domain decomposition to further increase the accuracy in the representation of the QoI. To illustrate the construction, we present numerical results for one-dimensional time varying linear and nonlinear diffusion equations with random space-dependent diffusion coefficients. Stochastic dimension reduction techniques proposed in the literature have mainly focused on quantifying the uncertainty in time independent and scalar QoI. To the best of our knowledge, this is the first time-dependent dimension reduction approach for time-dependent stochastic PDEs.
\end{abstract}

\begin{keyword}
Stochastic basis adaptation, time-dependent problems, domain decomposition 
\end{keyword}

\end{frontmatter}

\section{Introduction}

Uncertainty quantification for stochastic partial differential equations SPDEs with a large number of input random parameters is computationally expensive. Various stochastic dimension reduction methods e.g. basis adaptation \cite{tipireddy2014basis, Tipireddy2013, Tsilifis2016}, active subspace \cite{Constantine2014} and sliced inverse regression \cite{Li2016} have been developed to represent a scalar QoI in a low dimensional stochastic basis. In our prior work \cite{tipireddy2017basis, tipireddy2018stochastic} we have combined stochastic basis adaptation and spatial domain decomposition \cite{Chen2015, Toselli2005} to represent and compute the solution of a stochastic steady state PDE using a low dimensional basis. The low dimensional stochastic basis in each subdomain was obtained by Hilbert KL expansion \cite{Doostan2007}. 
Most physical systems are partially observed and evolve in time, and can be modeled by nonlinear time-dependent SPDEs. However, the existing stochastic dimension-reduction methods are limited to steady-state problems.  
Extending the dimension reduction methods to time-dependent problems is challenging because small errors in the low dimensional approximate solution introduced in the previous time steps  accumulate over time and cause divergence from the true solution. To address this challenge, we extend the basis adaptation and domain decomposition methods \cite{Tipireddy2013, tipireddy2014basis, tipireddy2017basis, tipireddy2018stochastic} to the case of time-dependent SPDEs. In the proposed approach,  the basis adaptation varies at every time step so that the low dimensional basis is adapted to the solution at that time step. 

As a first application to time-dependent problems, we apply it to a one-dimensional linear and a one-dimensional nonlinear (Richards) stochastic diffusion equations. Unlike the two- and three-dimensional computational domains where the interface between two subdomains is a line and a surface respectively, algorithms for domain decomposition in one-dimensional domain can be greatly simplified because the interface between two subdomains is a point. Hence, the unknown quantities at each interface will be two scalar quantities namely the solution and the flux. These scalar quantities can be obtained by solving a set of simple linear algebraic equations in the case of linear diffusion equation because of the superposition principle. However, for the nonlinear Richards equation we use an iterative approach to obtain these quantities. The proposed time-dependent basis adaptation can be extended to two and three-dimensional problems following the iterative domain decomposition algorithms such as Neumann-Neumann algorithm shown in \cite{tipireddy2018stochastic}.

The paper is organized as follows. In Section \ref{sec:spde}, we introduce the general type of SPDE that we are interested in, along with the necessary expansions for the source of stochasticity and the solution. In section \ref{sec:ba_time}, we introduce the time-dependent basis adaption method. In Section 4 we illustrate it through a linear time-dependent diffusion equation and a nonlinear time-dependent diffusion equation whose coefficients are modeled as log-normal random fields. Section \ref{sec:numerical} contains numerical results. Section \ref{sec:conclusions} contains conclusions and discussion of future work.  

\section{Time Dependent Stochastic PDEs} \label{sec:spde}
Let  $D$ be an open subset of $\mathbb{R}^n$ and $\Omega$ a
sample space. We want to find
find a random field, $u(x,\omega):D \times \Omega \rightarrow \mathbb{R}$ such that:
\begin{equation}\label{eq:spdeop}
 \mathcal{L}(x,u(t, x,\omega);a(x,\omega)) = f(x,\omega)  \;\; \rm{in}~D\times \Omega,
\end{equation}
subject to the boundary condition
\begin{equation}\label{eq:spdebc}
 \mathcal{B}(x,u;a(x,\omega)) = b(x,\omega)  \;\; \rm{on}~\partial D\times \Omega,
\end{equation}
where $\mathcal{L}$ is a differential operator and $\mathcal{B}$ is a boundary operator. We model the uncertainty in the stochastic PDE by treating the coefficient, $a(x,\omega)$ in the differential operator as a random field and propagate the uncertainty to the solution field, $u(x,\omega)$. To solve the stochastic PDE numerically, we discretize the random fields, $a(x,\omega)$ and $u(x,\omega),$ both in the spatial and stochastic domains. We model $a(x,\omega)$ with a lognormal random field, such that $a(x,\omega) = \exp[g(x,\omega)]$, where $g(x,\omega)$ is a Gaussian random field whose mean and covariance function are known. 

\subsection{Representation of the random coefficient and the solution}\label{sec:pce}

The Gaussian random field $g(x,\omega)$ can be approximated through a truncated KL expansion
\begin{equation}\label{eq:gkl}
g(x,\omega) \approx \tilde{g}(x,\xi(\omega)) = g_0(x) + \sum_{i=1}^M \sqrt{\lambda_i} g_i(x) \xi_i(\omega),
\end{equation}
where $g_0(x)$ is the mean of the random field $g(x,\omega)$ and $(\lambda_i, g_i(x))$ are the eigenvalues and eigenvectors obtained by solving the integral eigenvalue problem
\begin{equation}\label{eq:eig}
	\int_D C_g(x_1,x_2)  g_i(x_2)dx_2 = \lambda_i g_i(x_1).
\end{equation}
The eigenvalues are positive and non increasing, and the eigenfunctions $g_i(x)$ are orthonormal, 
\begin{equation}\label{eq:ortho}
	\int_D g_i(x) g_j(x)dx = \delta_{ij},
\end{equation}
where, $\delta_{ij}$ is the Kronecker delta. The random variables $\xi_i$ in Eq.~\ref{eq:gkl} are uncorrelated with zero mean. We model the random variables with uniform distribution. We further assume that they are independent. 

We approximate $a(x,\omega)$ and the solution field, $u(t, x,\omega)$ with truncated PC expansions with multi-variate orthogonal polynomials 

\begin{equation}\label{eq:pce_a}
a(x,\omega) \approx \tilde{a}(x,\boldsymbol{\xi}(\omega)) = a_0(x) + \sum_{i=1}^{N_{\xi}} a_i(x) \psi_i(\boldsymbol{\xi}),
\end{equation}

and

\begin{equation}\label{eq:pce_u}
u(t, x,\omega) \approx \tilde{u}(t, x,\boldsymbol{\xi}(\omega)) = u_0(t, x) + \sum_{i=1}^{N_{\xi}} u_i(t, x) \psi_i(\boldsymbol{\xi}),
\end{equation}
where, $\boldsymbol{\xi} = (\xi_1, \cdots, \xi_M)^T$, $u_0(x)$ is the mean of the solution filed, $u_i(x)$ are polynomial chaos coefficients and $\{\psi_i(\boldsymbol{\xi})\}$ are multi-variate Hermite polynomials. These polynomials are orthonormal with respect to the inner product defined by the expectation in the stochastic space, 

\begin{equation}\label{eq:innprod}
\langle \psi_i(\boldsymbol{\xi}), \psi_j(\boldsymbol{\xi}) \rangle \equiv \int_{\Omega} \psi_i(\boldsymbol{\xi}(\omega)) \psi_j(\boldsymbol{\xi}(\omega)) dP(\omega) = \delta_{ij}.
\end{equation}

\section{Time-dependent basis adaptation} \label{sec:ba_time}
In this section we describe the proposed approach that combines spatial domain decomposition and time-dependent stochastic basis adaptation.

\subsection{Basis adaptation at each time step} 
We will expand on our prior work on uncertainty quantification \cite{Tipireddy2013, tipireddy2014basis, tipireddy2017basis, tipireddy2018stochastic}, which combines  stochastic basis adaptation and spatial domain decomposition to the solution of time-dependent problems. For time-independent problems, the reduced stochastic basis is constructed only once and used to compute a low-dimensional representation of the stochastic solution. In the case of time-dependent problems, the stochastic basis constructed at one time instant is not enough for the accurate low-dimensional representation of the stochastic solution at a different time instant. If one insists on using the same basis for all times, then the accumulation of errors leads to divergence from the true solution. Hence, a new stochastic basis adapted to the solution at each time instant is required. 

We assume that the stochasticity of the solution of the stochastic PDE at every time instant can be accounted for through a, usually expensive, sparse grid collocation method. {\it For each time instant,  we can also use a low level sparse grid to obtain only the linear terms in a PC expansion of the solution in the entire spatial domain}. The linear PC expansion can be used to obtain the solution in each subdomain. Then, as is done in stochastic basis adaptation, we can use a Hilbert KL expansion in order to obtain a stochastic basis for each subdomain.


\subsection{Time-dependent basis adaptation in each subdomain}
We begin with the low level sparse grid solution, which allows us to compute the linear terms in the PC expansion of the solution in the {\it entire} spatial domain $D$ for the time duration of interest $T$ by solving~\eqref{eqn:1d_diff}, i.e., 
\begin{equation}\label{RK:eq:pce_ug2}
	 u_g(t, x,\boldsymbol{\xi}(\omega)) = u_0(t, x) + \sum_{i=1}^{d} u_i(t, x) \xi_i,
\end{equation}
where, $u_i(t, x) $ are the time-dependent PC coefficients and $\xi,\; i=1,\ldots,d$ are $d$ (sparse grid) collocation points. The computational cost of obtaining $u_g(t,x)$ is small compared to computing the full solution $u(t,x).$

We decompose the spatial domain $D$ into $K$ non-overlapping subdomains, $D^{(i)}, i=1,\ldots,K$ such that,
\begin{equation}\label{eq:dd}
	D = \cup_{i=1}^K D^{(i)}, \quad D^{(i)} \cap D^{(j)} = \emptyset
\end{equation}
Let the linear part of the solution in subdomain $D_s$ be  $u^s_g(t, x,\boldsymbol{\xi}(\omega)).$ We have
\begin{equation}\label{RK:eq:pce_ug_s}
	u^s_g(t, x,\boldsymbol{\xi}(\omega)) = u_g(t, x,\boldsymbol{\xi}(\omega)) \mathbb{I}_{D_s}(x), 
\end{equation}
where $\mathbb{I}_{D_s}(x)$ is the indicator function so that for any set $D_s,$ $\mathbb{I}_{D_s} = 1$ if $x \in D_s,$ and $\mathbb{I}_{D_s} = 0$ if $x \notin D_s$. 

We use the Hilbert space KL expansion \cite{Levy1999, Kirby1992, Silverman1996, Berkooz1993539, Christensen1999} of $u(t, x,\omega),$ to find the reduced basis $\boldsymbol{\eta}^s_t$ in each subdomain $D_s.$ To use the KL expansion, the solution should satisfy certain regularity and smoothness conditions, such as the solution being a subset of $L_2(\omega)$ \cite{Doostan2007}.

To find the time-dependent low dimensional basis in $D_s,$ we construct the covariance function of $u^s_g(t, x,\boldsymbol{\xi})$ in each subdomain $D_s$ as follows: 
\begin{equation}\label{RK:eq:cov_ug}
	 C^{s,t}_{u_g}(x_1,x_2) = \sum_{i=1}^{d} u_i(t, x_1) u_i(t, x_2), \quad x_1, x_2 \in {D_s}, t \in [0, T].
\end{equation}
The Hilbert space KL expansion of $u^s_g(t, x,\boldsymbol{\xi})$  \cite{Doostan2007}) in subdomain $D_s$ can be written as 
\begin{equation}\label{RK:eq:ugkl}
	u^s_g(t, x,\boldsymbol{\xi}(\omega)) = u^s_0(t, x) + \sum_{i=1}^{d} \sqrt{\mu^s_i(t)} \phi^s_i(t, x) \eta^s_i(t, \omega), \quad x \in {D_s}, t \in [0,T],
\end{equation}
where $(\mu^s_i(t), \phi^s_i(t, x))$ are time-dependent eigenvalue and eigenfunction in the Hilbert KL expansion obtained by solving the eigenvalue problem at a time step $t$
\begin{equation}\label{eq:cov_cg}
	 \int_{D_s} C^{s,t}_{u_g}(x_1,x_2) \phi^s_i(t, x_1)dx_1 = \mu^s_i(t) \phi^s_i(t, x_2), \quad i = 1, 2, \ldots, d. 
\end{equation}
The random variable $\eta^s_i$ can be written as
\begin{align}\label{eq:pce_eta}
	 \eta^s_i(t) &= \frac{1}{\sqrt{\mu^s_i(t)}}\int_{D_s} \left(u^s_g(t, x,\boldsymbol{\xi}) - u^s_0(t, x) \right) \phi^s_i(t, x)dx, \nonumber \\
	 		&= \frac{1}{\sqrt{\mu^s_i(t)}}\int_{D_s} \left(u_0(t, x) + \sum_{j=1}^{d} u_j(t, x) \xi_j - u^s_0(t, x) \right) \phi^s_i(t, x)dx. \quad i = 1, 2, \ldots, d.
\end{align}
By construction, $u_0(t, x) = u^s_0(t, x)$ for $x\in D_s,$  and we find 
 \begin{align}\label{eq:eta_xi}
	 \eta^s_i(t) &= \frac{1}{\sqrt{\mu^s_i(t)}}\int_{D_s} \left ( \sum_{j=1}^{d} u_j(t,x) \xi_j \right) \phi^s_i(t,x)dx, \quad x \in D_s, t \in [0, T], i = 1, 2, \ldots, d \nonumber \\
	 		&= \sum_{j=1}^{d} \left ( \frac{1}{\sqrt{\mu^s_i(t)}}\int_{D_s}  u_j(t,x) \phi^s_i(t,x)dx  \right)\xi_j, \quad x \in D_s, t \in [0,T], i = 1, 2, \ldots, d \nonumber \\
	 		&= \sum_{j=1}^{d} a^s_{ij}(t) \xi_j,
\end{align}
where $a^s_{ij(t)} =  \frac{1}{\sqrt{\mu^s_i(t)}}\int_{D_s}  u_j(t,x) \phi^s_i(t,x)dx, i,j = 1, \ldots, d.$ The quantities $a^s_{ij}(t)$ provide a time-dependent linear map between the random variables $\eta^s_i(t)$ and the original random variables $\xi_j,$ for $j=1,\ldots,d.$ In each subdomain, and at each time instant $t$, we can use the reduced dimensional basis $\{\eta^s_i\}$ to solve the stochastic PDE in question using sparse grid collocation which is a non-intrusive method. 



\section{Domain decomposition and basis adaptation for one-dimensional diffusion equations}
In this section, we present two test problems namely, a linear stochastic diffusion equation and nonlinear stochastic diffusion (Richards) equation for demonstrating the proposed time dependent basis adaptation and domain decomposition methods. 

\subsection{Linear diffusion equation}\label{sec:dd}
Here, we use a linear 1D time-dependent stochastic diffusion equation to illustrate the construction. We consider 
\begin{align} \label{eqn:1d_diff}
	& \frac{\partial u(t, x, \omega)}{\partial t} - \frac{\partial}{\partial x}\left [ a(x, \omega) \frac{\partial}{\partial x} u(t, x, \omega)   \right ] = f, \quad \quad 0<x<L 
\end{align}
with boundary conditions
\begin{align} \label{RK:eqn:rich_boundary}
    u(t,0) = u_0 \quad \text{and}  \quad u(t, L)=u_L,
\end{align}
and initial condition $u(0,x) = u_0(x).$  Due to the linearity of the problem we will use the superposition principle to compute the solution after domain decomposition. 

For each of the non-overlapping subdomains, $D^{(i)}, i=1,\ldots,K$ we solve
%
\begin{align} \label{eqn:1d_diff_dd}
	& \frac{\partial u^{(i)}(t, x, \omega)}{\partial t} - \frac{\partial}{\partial x}\left [ a(x, \boldsymbol{\xi}) \frac{\partial}{\partial x} u^{(i)}(t, x, \omega)   \right ] = f, \quad \quad x \in D_s
\end{align}
with interface and boundary conditions
\begin{align} \label{RK:eqn:rich_nl}
    \mathcal{C}^{(ij)}(x,u^{(i)},u^{(j)};a(x,\omega)) = \beta(x,\omega)  \;\; \rm{on}~\partial D^{(i)} \cap \partial D^{(j)} \times \Omega, \\ 
    u^{(i)}(t,0) = u_0 \quad \text{and}  \quad u^{(i)}(t, L)=u_L, \rm{on}~\partial D^{(i)} \cap  \partial D\times \Omega.
\end{align}
Since this is a one-dimensional linear equation, the domain decomposition method can be simplified because the solution at the interface is a scalar quantity. The solution in each subdomain can be obtained using $0$ and $1$ boundary conditions at the interfaces to obtain the solutions. Then, we can use the superposition principle and apply interface conditions to get the unknown solution at the interfaces.


\subsection{Nonlinear diffusion equation} \label{sec:1d_richards}
Here, we consider the nonlinear diffusion (Richards) equation with the van Genuchten model \cite{pan1995transformed}:

\begin{align} \label{eqn:1d_richards}
	& \frac{\partial \psi(t, x, \omega)}{\partial t} - \frac{\partial}{\partial x}\left [ K(x, \omega) \left(\frac{\partial}{\partial x} \psi(t, x, \omega) +1 \right)  \right ] = 0, \quad \quad 0<x<L 
\end{align}
with boundary conditions
\begin{align} \label{RK:eqn:rich_boundary2}
    \psi(t,0) = \psi_0 \quad \text{and}  \quad \psi(t, L)=\psi_L,
\end{align}
and initial condition $\psi(0,x)$, where $\psi$ is the pressure head. The water retention $S_e$ and hydraulic conductivity $K$ are related to the pressure head $\psi$ \cite{tipireddy2018stochastic} through
\begin{equation} \label{eq:vg_Se}
    S_e = \frac{\theta - \theta_r}{\theta_s - \theta_r} = \left[\frac{1}{1+(\alpha_{vg}|\psi|)^n} \right]^m,
\end{equation}
and
\begin{equation} \label{eq:vg_K}
    K = K_s\sqrt{S_e}[1-(1-S_e^{1/m})^m]^2,
\end{equation}
where $m=1-1/n;$ $\theta_r$ and $\theta_s$ are the residual and saturated water contents respectively; $\alpha_{vg}$  and $n$ are the van Genuchten model parameters. In this equation, we model the uncertain hydraulic conductivity as a log-normal random field $K_s(x, \boldsymbol{\xi}).$ 
We note that since the Richards equation is nonlinear, obtaining the value at the interface between subdomains cannot be achieved through the superposition principle. Instead, we employ an iterative  algorithm where we start with an initial guess for the interface solution and update it at each iteration in order to satisfy the governing equation.

%

\section{Numerical results}\label{sec:numerical}
In this section, we present numerical solutions of a linear stochastic diffusion equation and a nonlinear stochastic diffusion equation with the domain decomposition and time-dependent basis adaptation method and provide comparison with the solutions obtained with the time-independent basis adaptation method. For the nonlinear diffusion equation, we also present solutions obtained without domain decomposition in order to assess the importance of domain decomposition. 

\subsection{Linear diffusion equation}
Consider \eqref{eqn:1d_diff} on the spatial domain $x\in(0,1)$ with the boundary conditions
\begin{align} \label{RK:eqn:rich_nl_boundary}
    u(t,0) = 2.0 \quad \text{and}  \quad u(t, 1)= 1.0,
\end{align}
and initial conditions
\begin{align} \label{RK:eqn:rich_nl_initial}
    u(0,x)= 
\begin{cases}
    2.0,& \text{if } 0.375 \leq x\leq 0.625\\
    1.0, & \text{otherwise}.
\end{cases}
\end{align}

\begin{figure}[!ht]
\begin{center}
{\includegraphics[scale=.65]{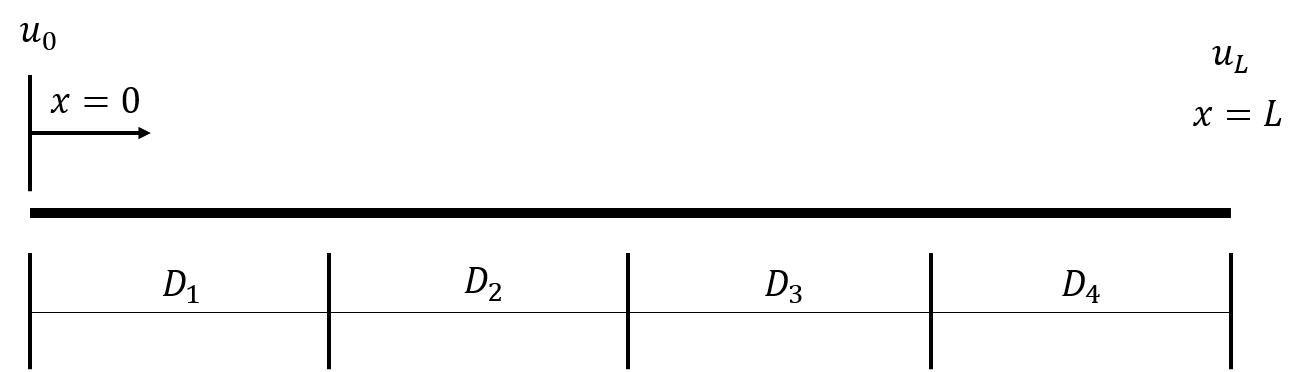}}
 \caption{Spatial domain decomposed into four subdomains} \label{fig:one_d_DD}
\end{center}
\end{figure} 
We model the random coefficient $a(x,\omega) = \exp[g(x,\omega)]$ using a log-normal distribution with mean $a_0 = 0.5$ and standard deviation $\sigma_a = 0.25.$ The correlation length of $a(x,\omega)$ is set to $0.2$. For the reference solution, we approximate the random coefficient $a(x,\omega)$ using KL expansion with $d=10$ random variables, i.e. $\boldsymbol{\xi}=(\xi_1,\ldots,\xi_{10}).$ To obtain the reference solution, we use the sparse-grid collocation method with sparse-grid level $5.$ The total number of deterministic simulations required to compute the reference solution is $8761.$ For the time-dependent basis adaptation, we decompose the spatial subdomain into four subdomains (see Fig. \ref{fig:one_d_DD}) and in each subdomain we approximate the solution using a stochastic basis of dimension $r=3.$ To obtain the time-dependent basis, we need to solve the time dependent SPDE at a few collocation points using a low sparse-grid level in full dimension. Here, we use the sparse-grid level 2 and $r=10$, which contain
$21$ collocation points. We use the solution at these collocation points to compute the linear transformation matrix $A$ needed to obtain the reduced basis $\boldsymbol{\eta} = A \boldsymbol{\xi}.$

\begin{figure}[ht!]
    \centering
    \begin{subfigure}[t]{0.32\textwidth}
        \centering
        \includegraphics[scale=.24]{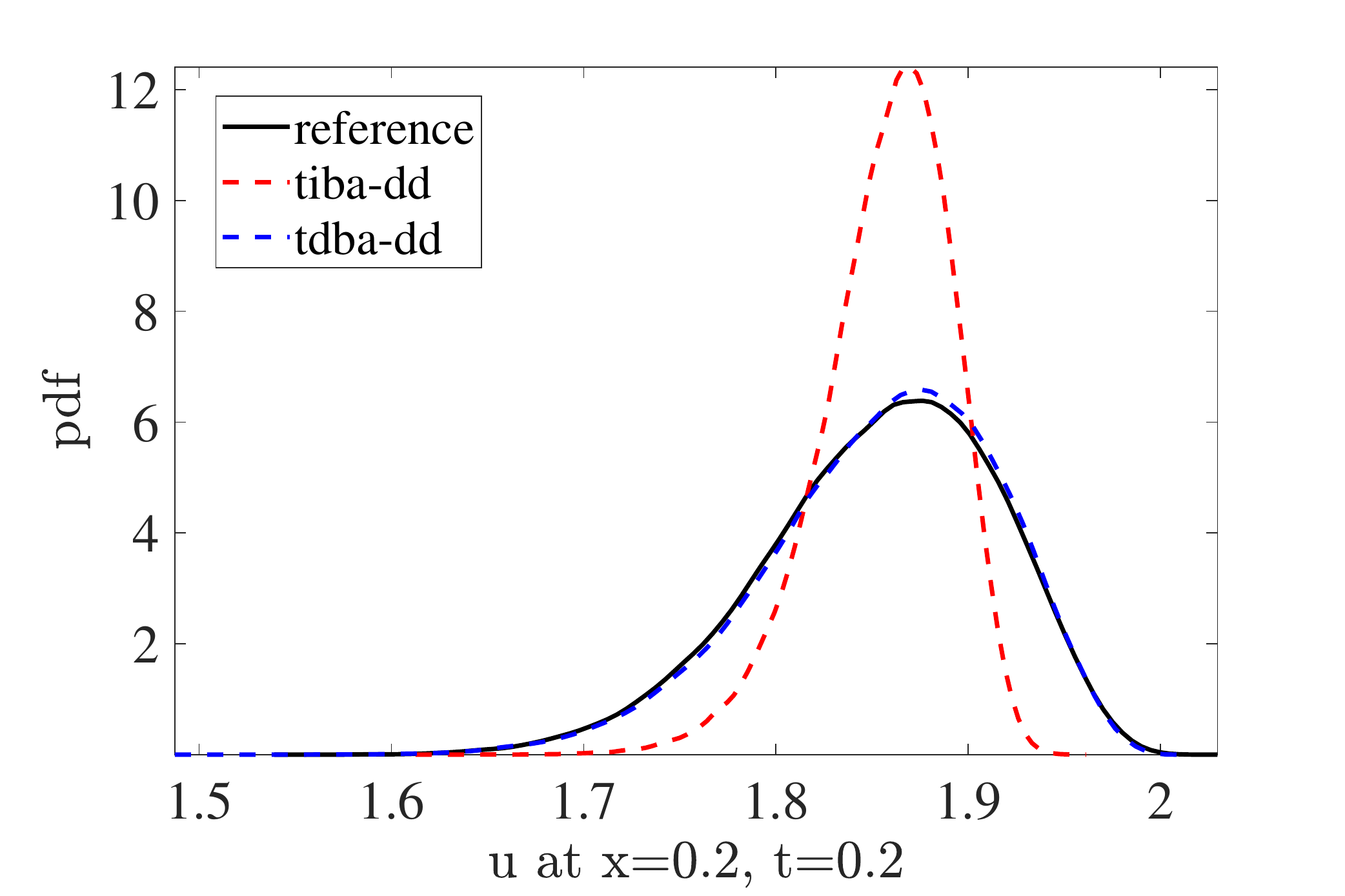}
        \caption{} \label{fig:pdf_u_at_t02_x02_dd_xid10_p3_etad3_p3}
    \end{subfigure}        
    \begin{subfigure}[t]{0.32\textwidth}
        \centering
        \includegraphics[scale=.24]{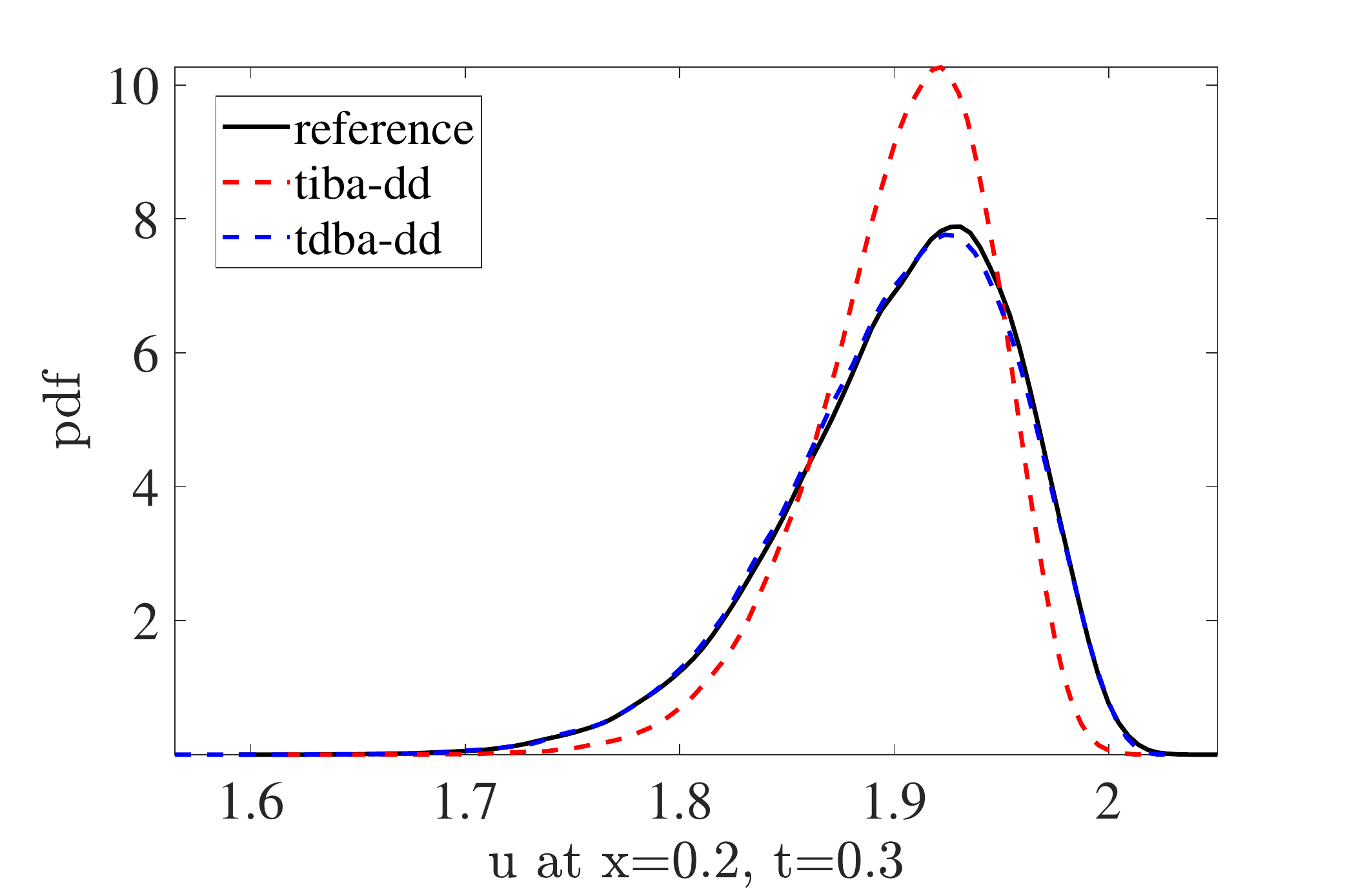}
        \caption{} \label{fig:pdf_u_at_t03_x02_dd_xid10_p3_etad3_p3}
    \end{subfigure} 
    \begin{subfigure}[t]{0.32\textwidth}
        \centering
        \includegraphics[scale=.24]{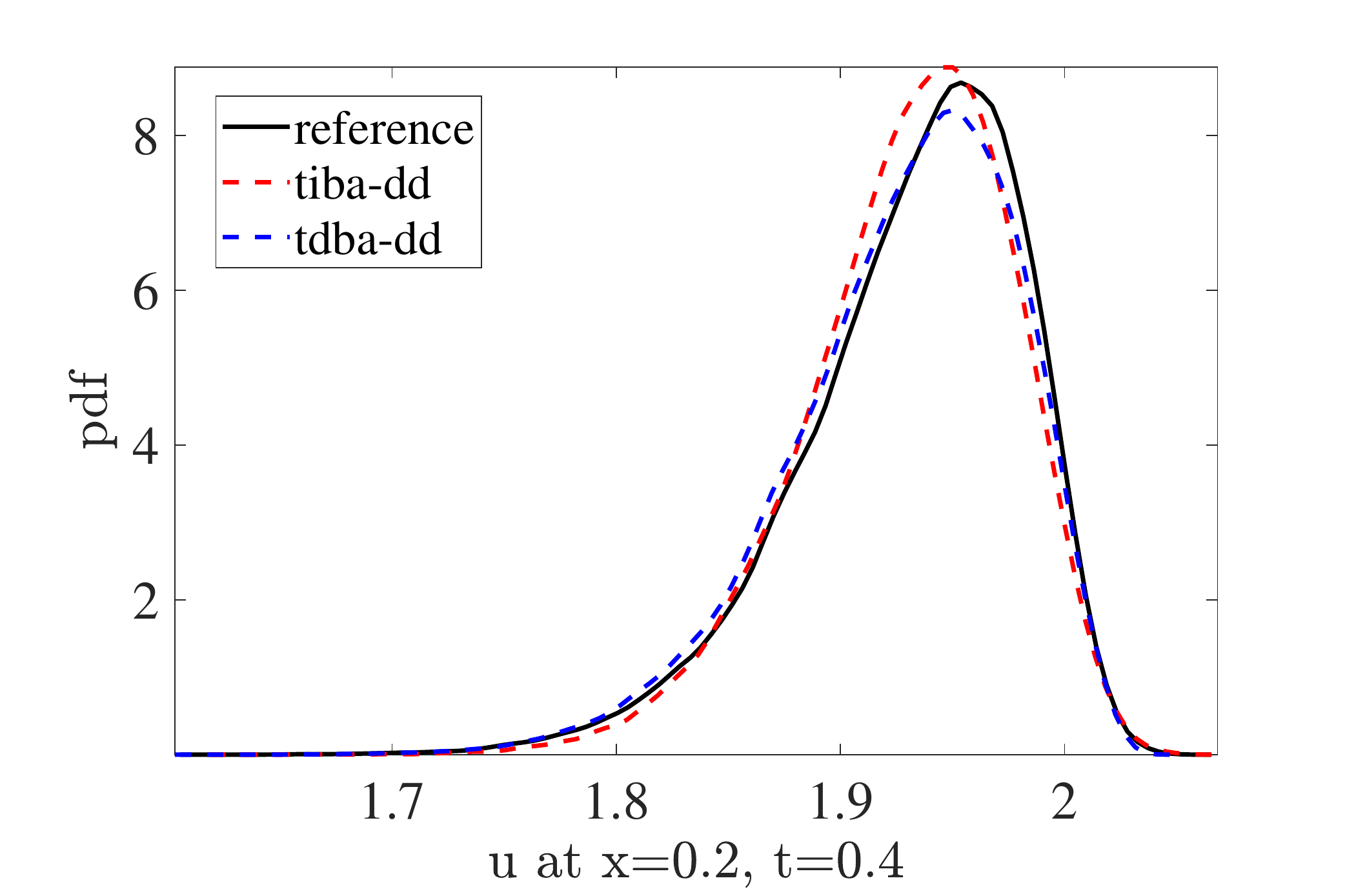}
        \caption{} \label{fig:pdf_u_at_t04_x02_dd_xid10_p3_etad3_p3}
    \end{subfigure}
    \begin{subfigure}[t]{0.32\textwidth}
        \centering
        \includegraphics[scale=.24]{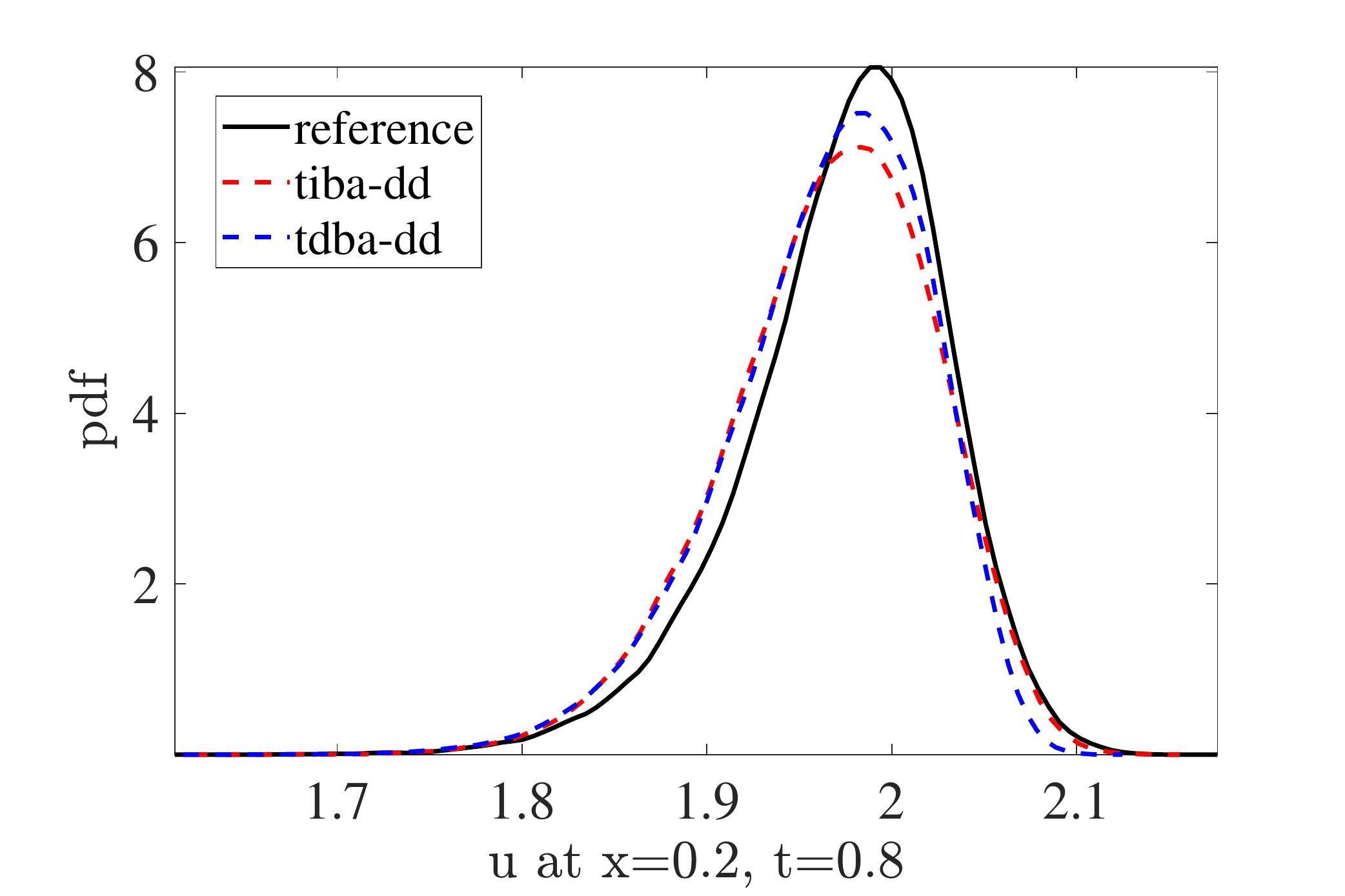}
        \caption{} \label{fig:pdf_u_at_t08_x02_dd_xid10_p3_etad3_p3}
    \end{subfigure}        
    \begin{subfigure}[t]{0.32\textwidth}
        \centering
        \includegraphics[scale=.24]{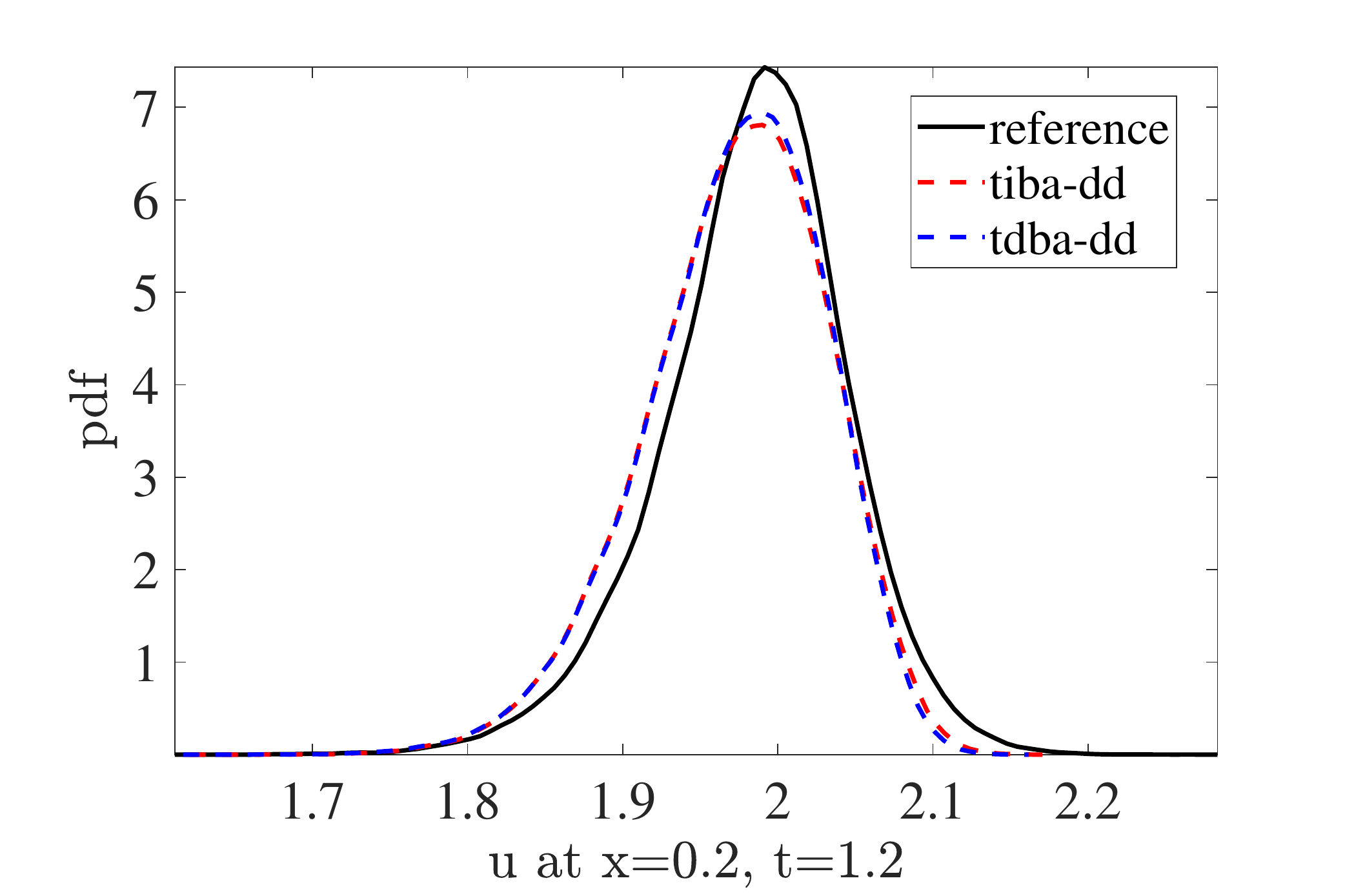}
        \caption{} \label{fig:pdf_u_at_t012_x02_dd_xid10_p3_etad3_p}
    \end{subfigure} 
    \begin{subfigure}[t]{0.32\textwidth}
        \centering
        \includegraphics[scale=.24]{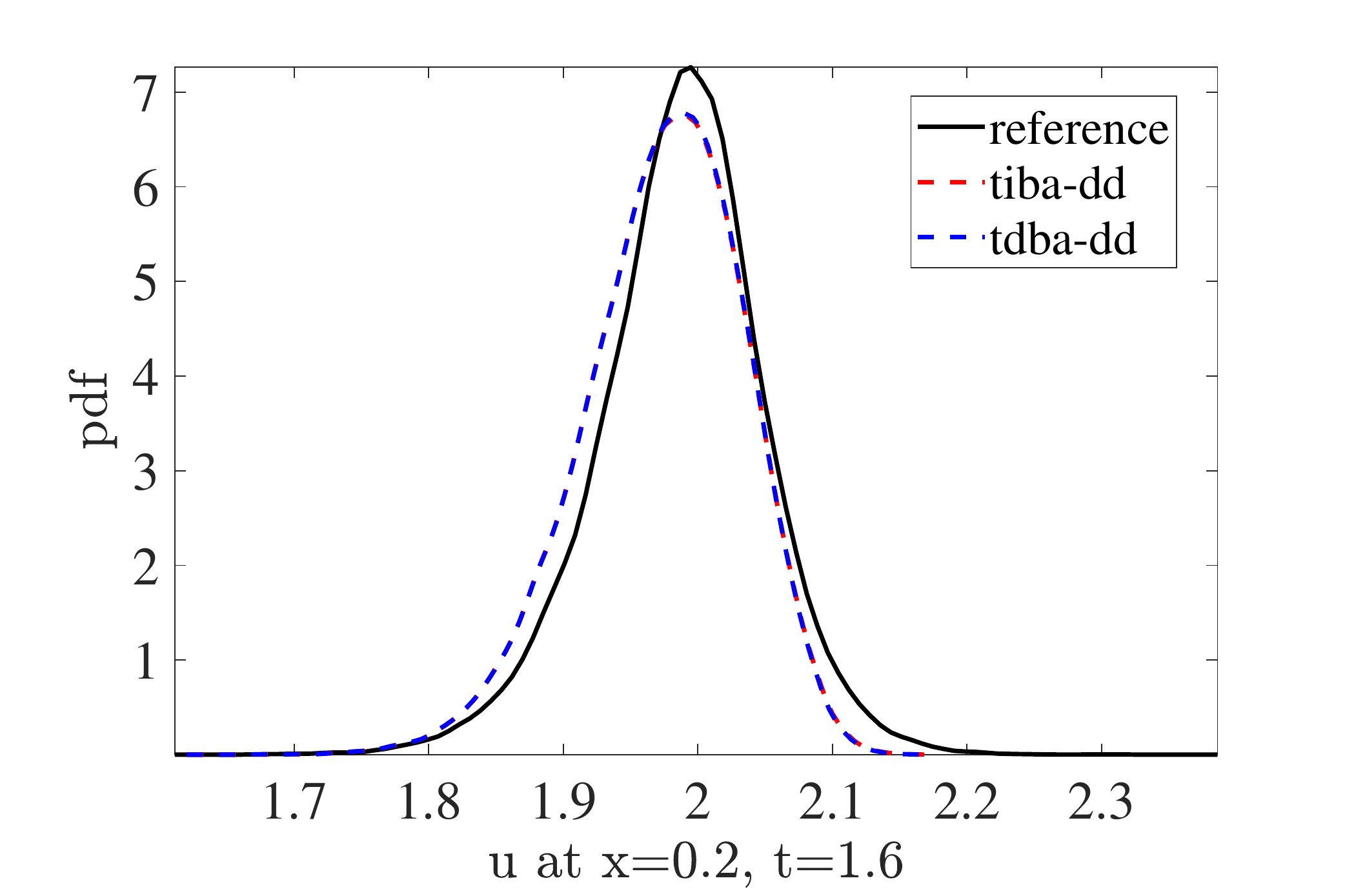}
        \caption{} \label{fig:pdf_u_at_t016_x02_dd_xid10_p3_etad3_p3}
    \end{subfigure}
\caption{Linear diffusion equation. Probability density function of $u$ at $x=0.7$ and a). $t=0.2$, b). $t=0.3$, c). $t=0.4$, d). $t=0.8$, e). $t=1.2$ and f). $t=1.6.$ The reference solution (solid black) is computed with dimension $d=10$ and the reduced dimensional solutions (dashed red and blue) are computed with $r=3$.} \label{fig:pdf_u_at_t_02_08_1_6_x07_dd_xid10_p3_etad3_p3}
\end{figure}

Fig. \ref{fig:pdf_u_at_t_02_08_1_6_x07_dd_xid10_p3_etad3_p3} shows the probability density function (pdf) of the solution at $x=0.7$ and time $t=0.2, 0.8$ and $1.6.$ We can see that for early times, the estimated pdf of the solution using time-dependent basis adaptation and domain decomposition (dashed blue line) matches well with the pdf of the reference solution (solid black line) whereas the pdf of the solution obtained with fixed basis adaptation and domain decomposition deviates (dashed red line). For longer times, the pdf estimates of the time-dependent and fixed basis adaptation converge. This is to be expected since the pdf of the diffusion equation converges to a steady state. In particular, while the choice of the basis adaptation is important for short times, eventually the contracting nature of the linear diffusion evolution operator makes both time-dependent and fixed basis adaptations to perform similarly. However, we note that the relatively inexpensive determination of the time-dependent basis adaptation allows us to follow the evolution of the pdf accurately for both short and long times. This can be seen clearly in Fig.2 where the accuracy of the time-dependent basis adaptation remains practically the same (note the change in scale in the vertical axis as time progresses).


\begin{figure}[ht!]
    \centering
    \begin{subfigure}[t]{0.48\textwidth}
        \centering
        \includegraphics[scale=.24]{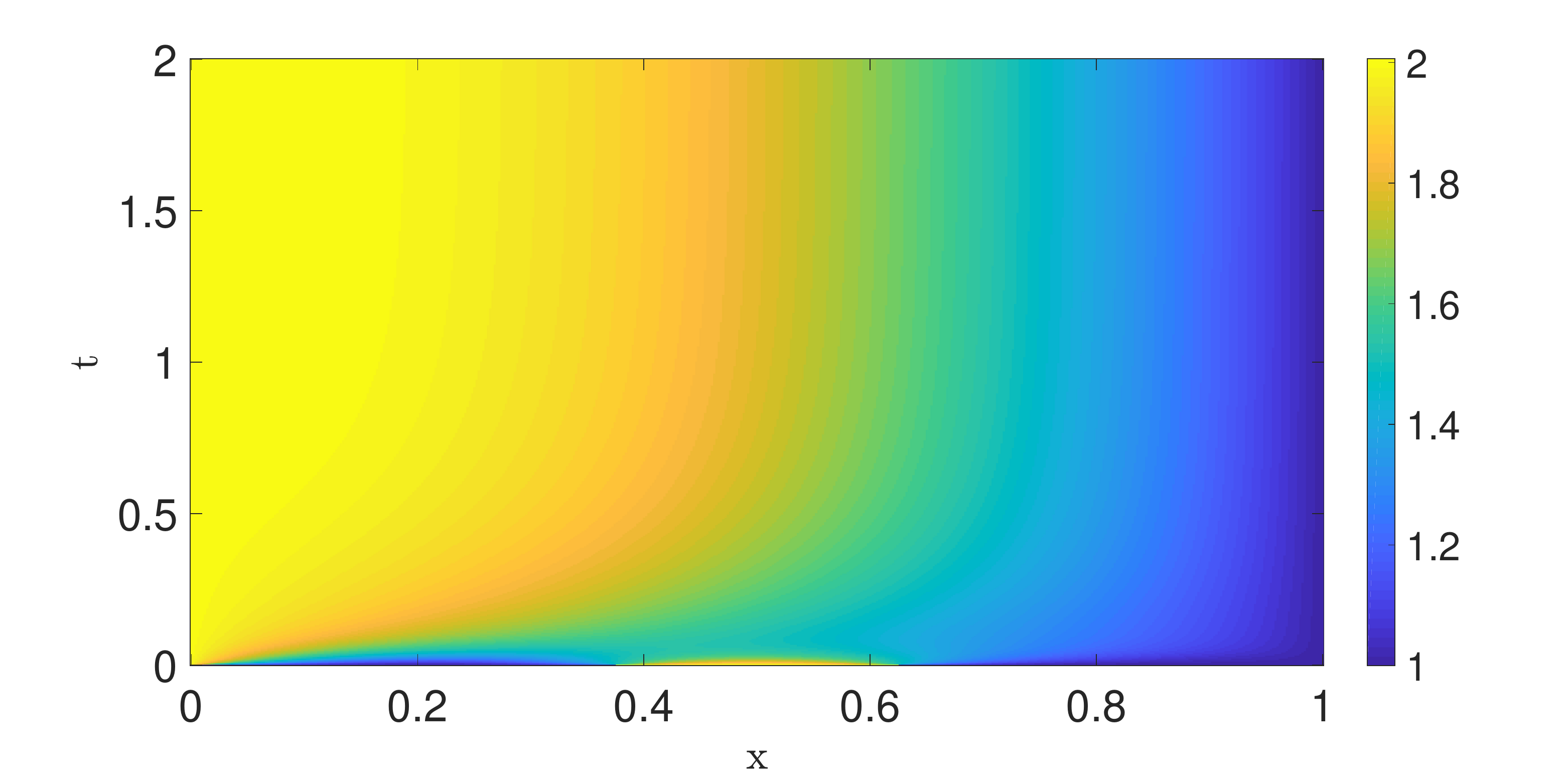}
        \caption{Mean} \label{fig:Uxt_mean_xi_10}
    \end{subfigure}        
    \begin{subfigure}[t]{0.48\textwidth}
        \centering
        \includegraphics[scale=.24]{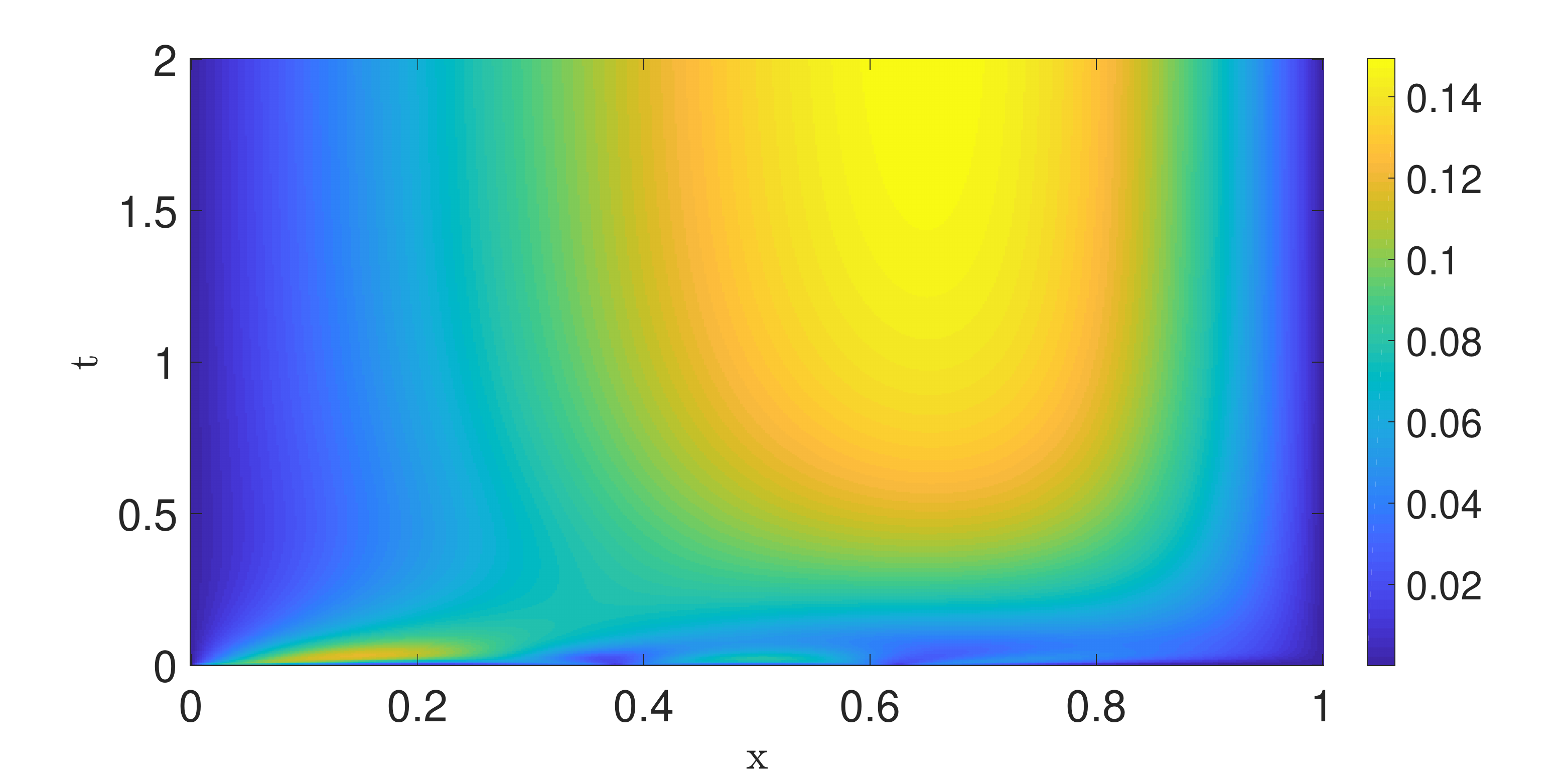}
        \caption{Standard deviation} \label{fig:Uxt_sdev_xi_10}
    \end{subfigure} 
    \begin{subfigure}[t]{0.48\textwidth}
        \centering
        \includegraphics[scale=.24]{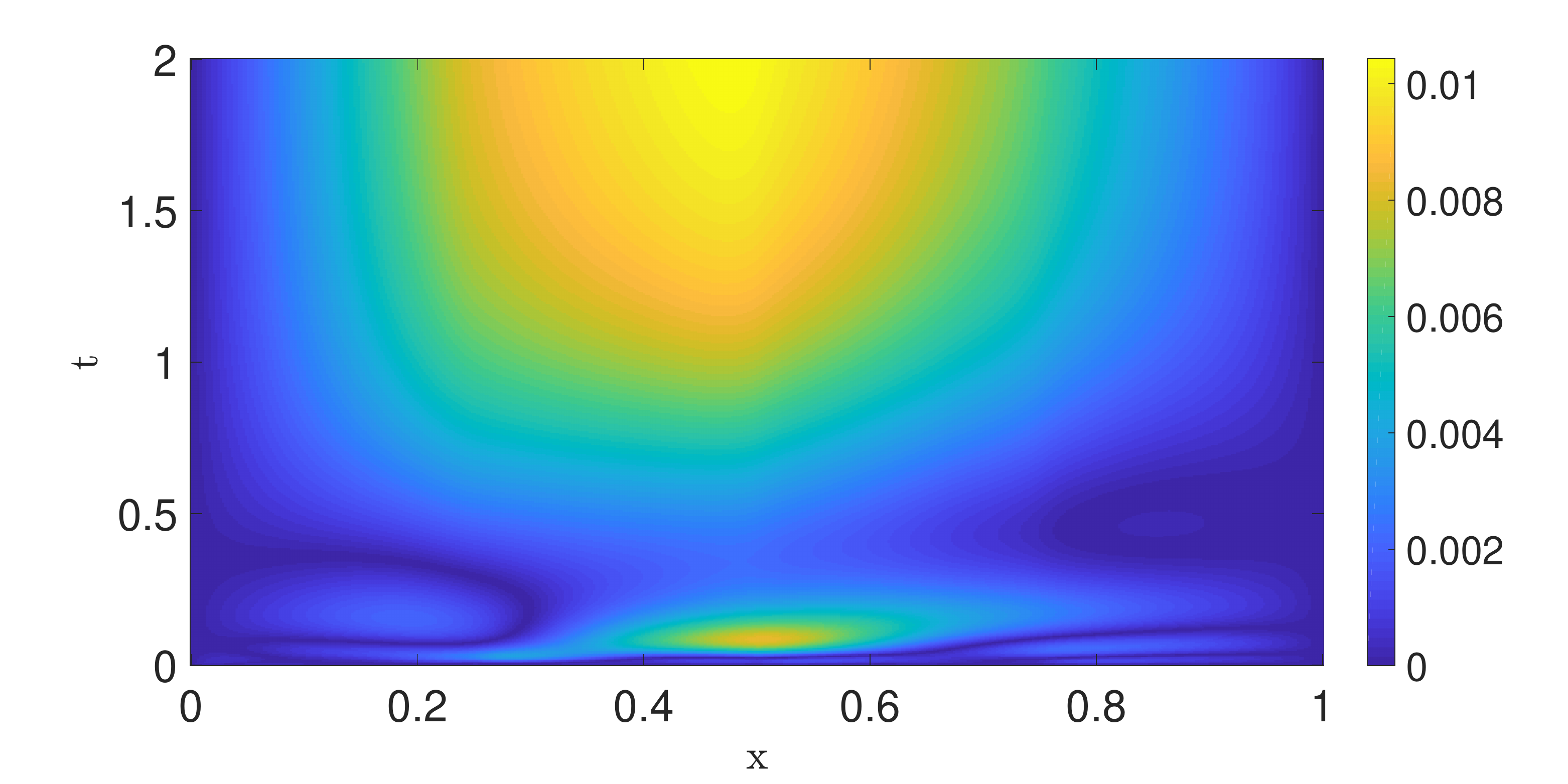}
        \caption{Abs. point. error in mean} \label{fig:Uxt_mean_rel_error_eta_3_tind}
    \end{subfigure}        
    \begin{subfigure}[t]{0.48\textwidth}
        \centering
        \includegraphics[scale=.24]{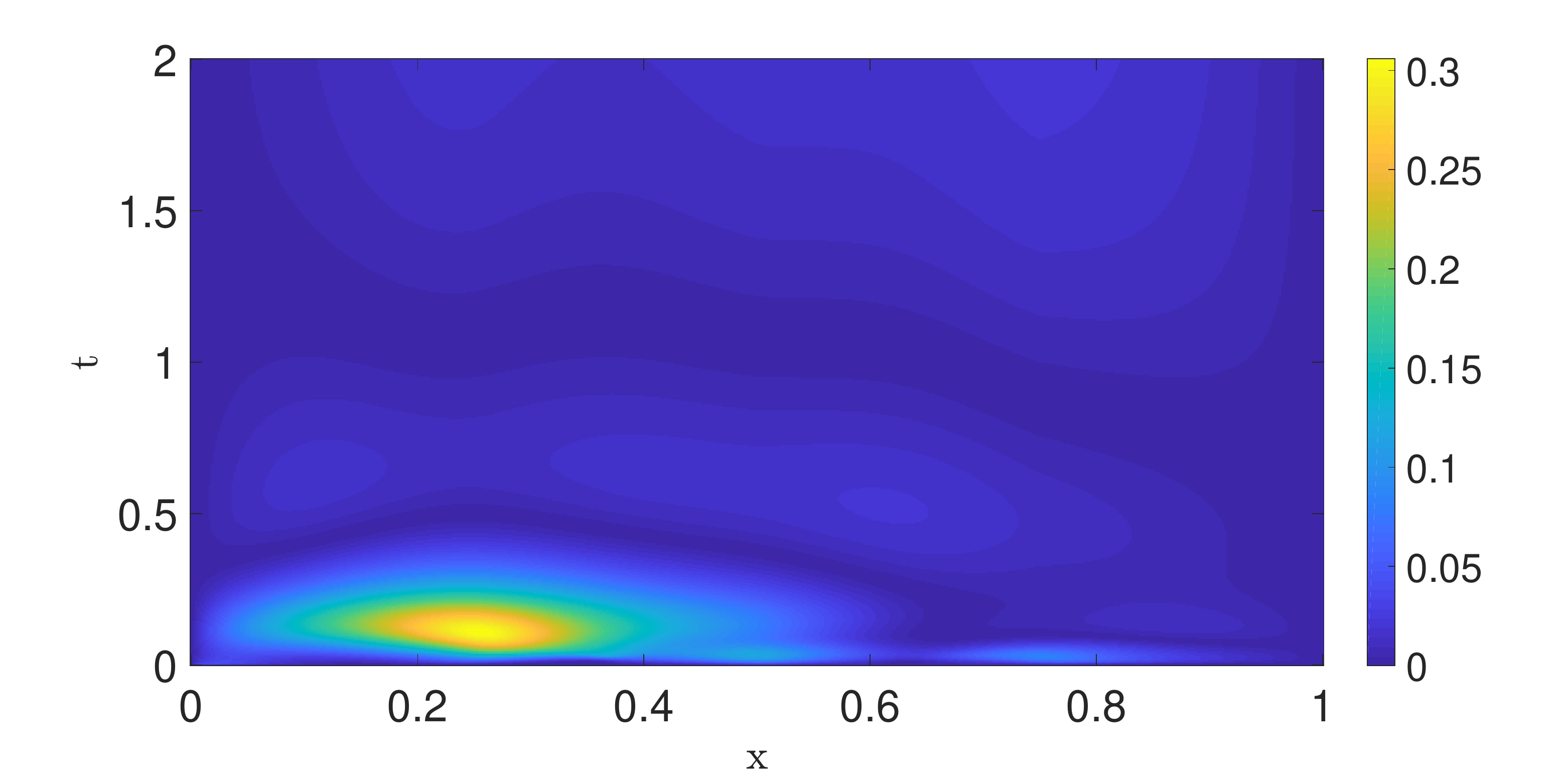}
        \caption{Abs. point error in standard deviation} \label{fig:Uxt_sdev_rel_error_eta_3_tind}
    \end{subfigure} 
    \begin{subfigure}[t]{0.48\textwidth}
        \centering
        \includegraphics[scale=.24]{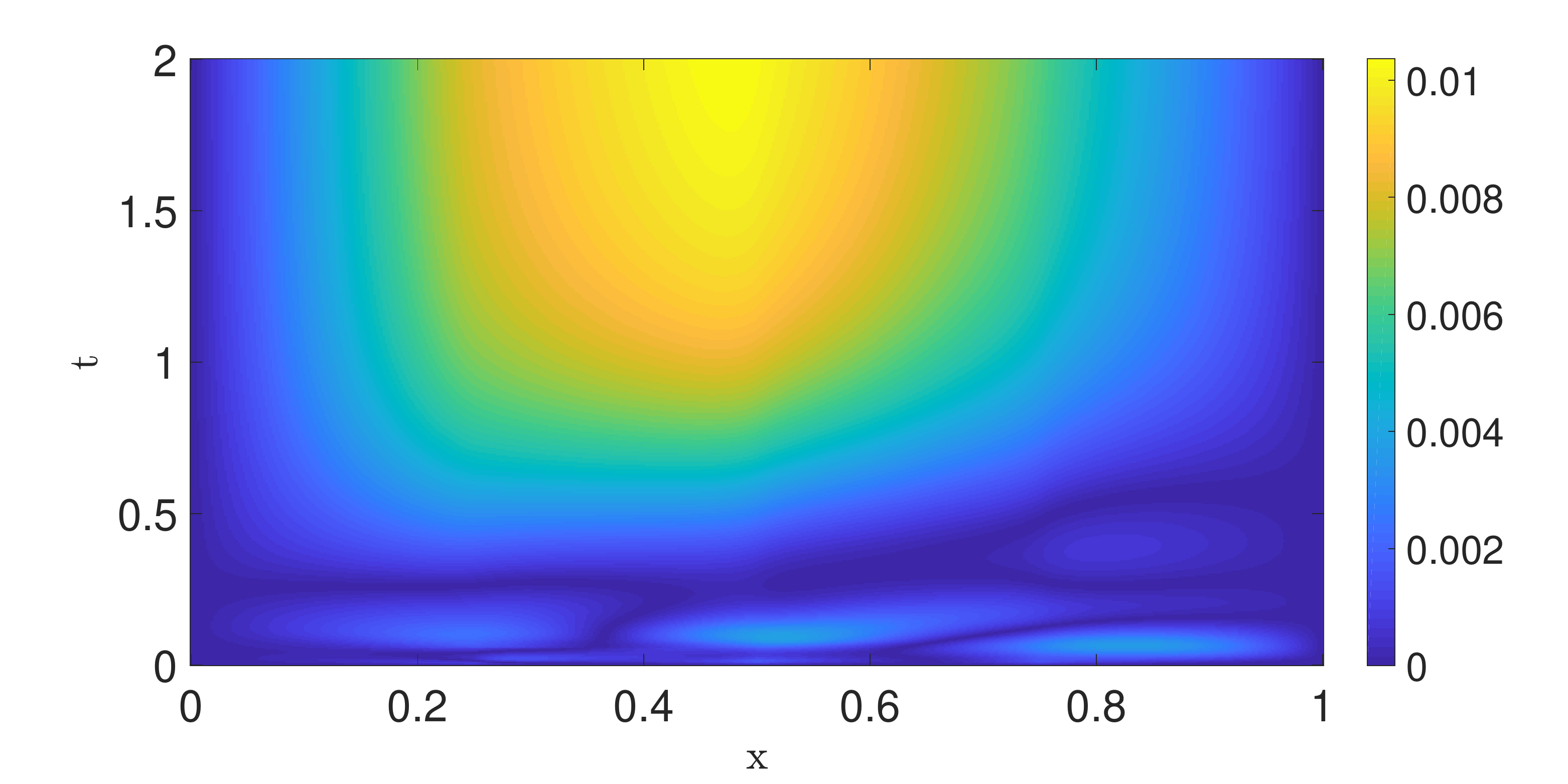}
        \caption{Abs. point error in mean} \label{fig:Uxt_mean_rel_error_eta_3_td}
    \end{subfigure}        
    \begin{subfigure}[t]{0.48\textwidth}
        \centering
        \includegraphics[scale=.24]{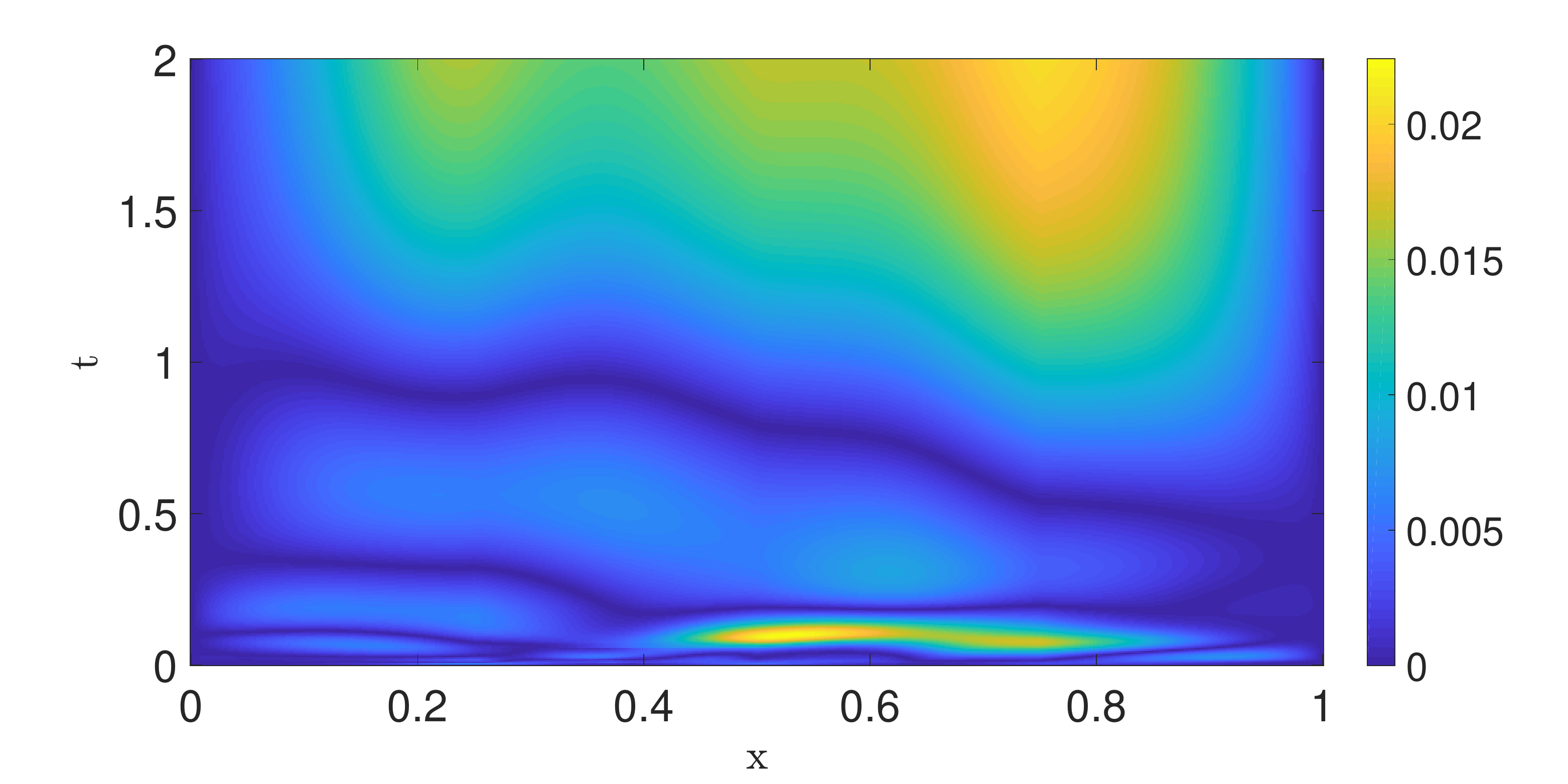}
        \caption{Abs. point error in standard deviation} \label{fig:Uxt_sdev_rel_error_eta_3_td}
    \end{subfigure} 
\caption{a) mean of reference solution ($d=10$) and b) standard deviation of reference solution ($d=10$) c) error in mean of low dimensional solution ($r=3$) with fixed basis adaptation and d) error in standard deviation of low dimensional solution ($r=3$) with fixed basis adaptation e) error in mean of low dimensional solution ($r=3$) with time-dependent basis adaptation and f) error in standard deviation of low dimensional solution ($r=3$) with time-dependent basis adaptation as a function of $x$ (horizontal axis) and $t$ (vertical axis).} 
\end{figure}

Figs. \ref{fig:Uxt_mean_xi_10} and \ref{fig:Uxt_sdev_xi_10} show the color map of the mean and standard deviation of the reference solution as a function of $x$ (horizontal axis) and $t$ (vertical axis). Figs. \ref{fig:Uxt_mean_rel_error_eta_3_tind} and \ref{fig:Uxt_sdev_rel_error_eta_3_tind} show the absolute point error in the mean and standard deviation of the solution using fixed basis adaptation. 
Figs. \ref{fig:Uxt_mean_rel_error_eta_3_td} and \ref{fig:Uxt_sdev_rel_error_eta_3_td} show the absolute point error of the mean and standard deviation of the solution using time-dependent basis adaptation compared to the reference solution. 

While both fixed and time-dependent basis adaptation methods result in similar relative error for the mean, time-dependent basis adaptation results in significantly lower relative error for the standard deviation. This is important because a smaller error in the standard deviation estimate improves the confidence in the prediction of the mean. 

\subsection{Nonlinear diffusion equation}
In Eqs (\ref{eqn:1d_richards})-(\ref{eq:vg_K}), we set $L=10$, $\psi_0 = 0$ and $\psi_L = -0.35$ and $\psi(0,x) = \frac{\psi_L-\psi_0}{x_L-x_0}x$. We model the coefficient $K_s(x,\boldsymbol{\xi})$ as a log-normal random field, $K_s(x,\boldsymbol{\xi})=\exp[g_s(x,\boldsymbol{\xi})]$, where $g_s$ is a Gaussian random field with the mean $\mu_k=5.0$, standard deviation $\sigma_k=0.5$, and the exponential correlation function with the correlation length $L/4$.
 The corresponding coefficient of variation is $\sigma_k/\mu_k=0.1$ indicating a weak heterogeneity. In this section, we demonstrate that even for such a small coefficient of variation, a significant error in the standard basis adaptation method (with constant basis) accumulates over time, while a solution obtained using the time-dependent basis matches well with the reference solution for the considered time interval.

\begin{figure}[ht!]
    \centering
    \begin{subfigure}[t]{0.32\textwidth}
        \centering
        \includegraphics[scale=.25]{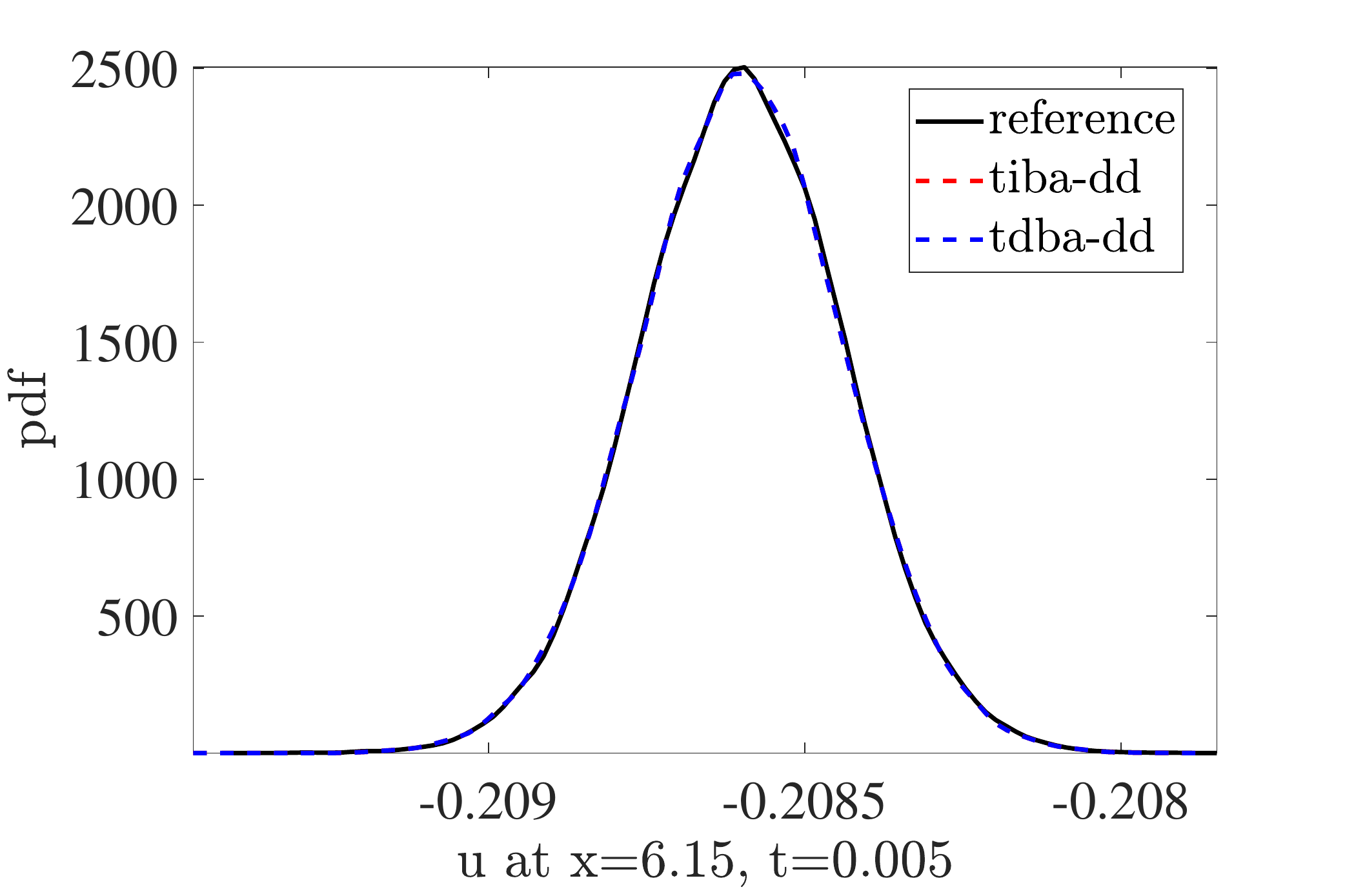}
        \caption{} \label{fig:pdf_u_at_t2_x62_dd_xid10_p3_etad5}
    \end{subfigure}        
    \begin{subfigure}[t]{0.32\textwidth}
        \centering
        \includegraphics[scale=.25]{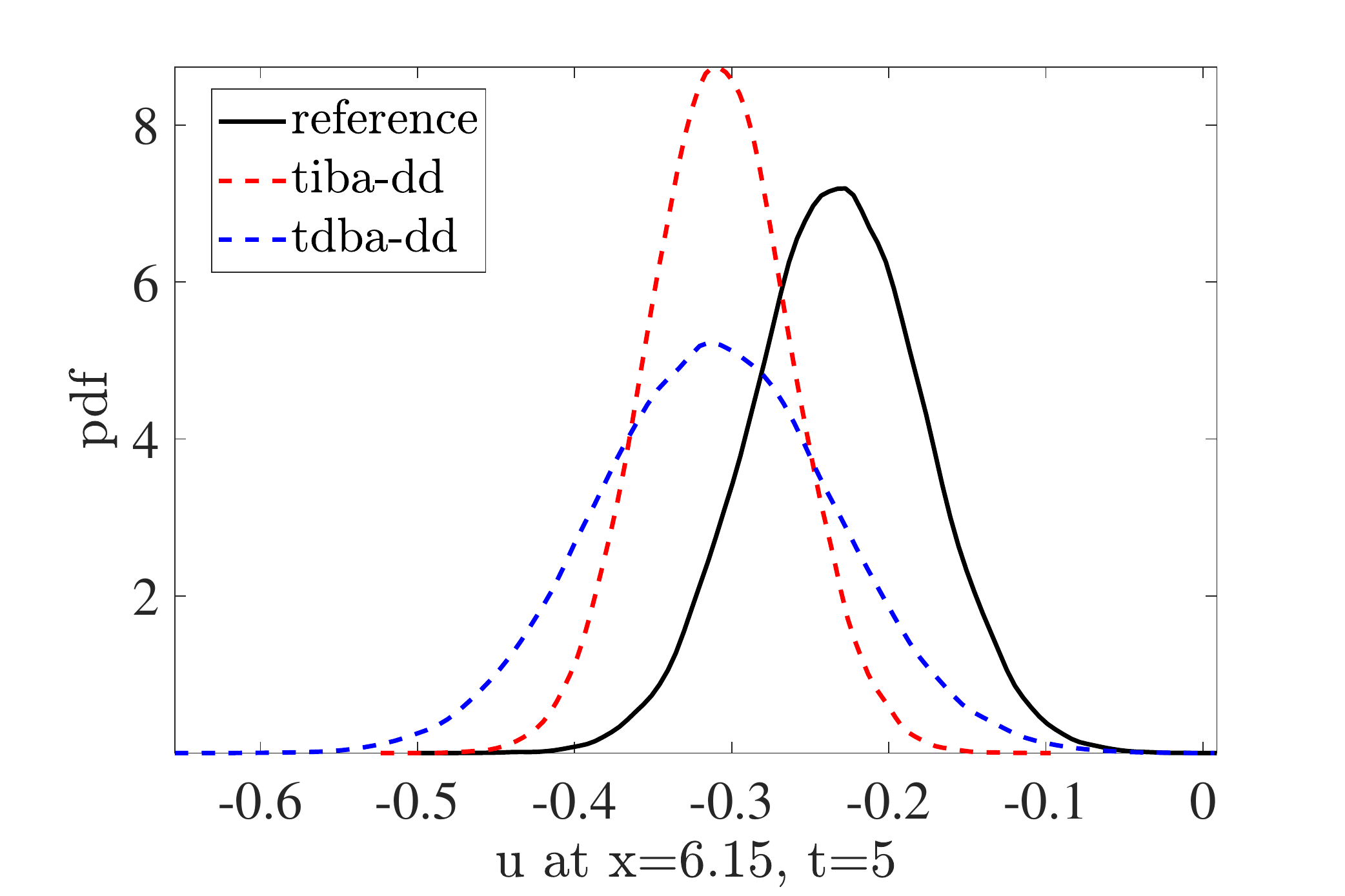}
        \caption{} \label{fig:pdf_u_at_t1001_x62_dd_xid10_p3_etad5}
    \end{subfigure} 
    \begin{subfigure}[t]{0.32\textwidth}
        \centering
        \includegraphics[scale=.25]{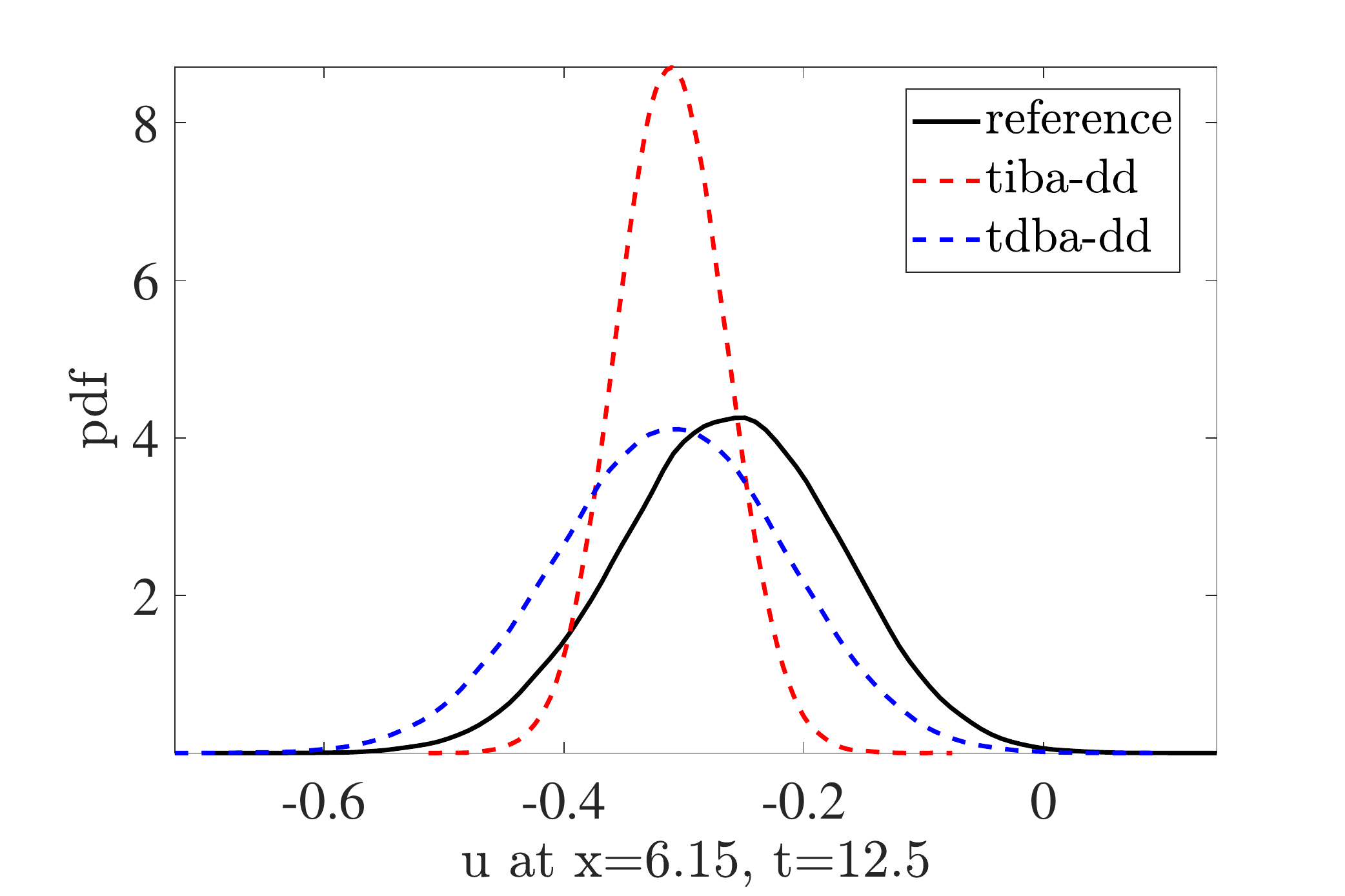}
        \caption{} \label{fig:pdf_u_at_t2501_x62_dd_xid10_p3_etad5}
    \end{subfigure}
    \begin{subfigure}[t]{0.32\textwidth}
        \centering
        \includegraphics[scale=.25]{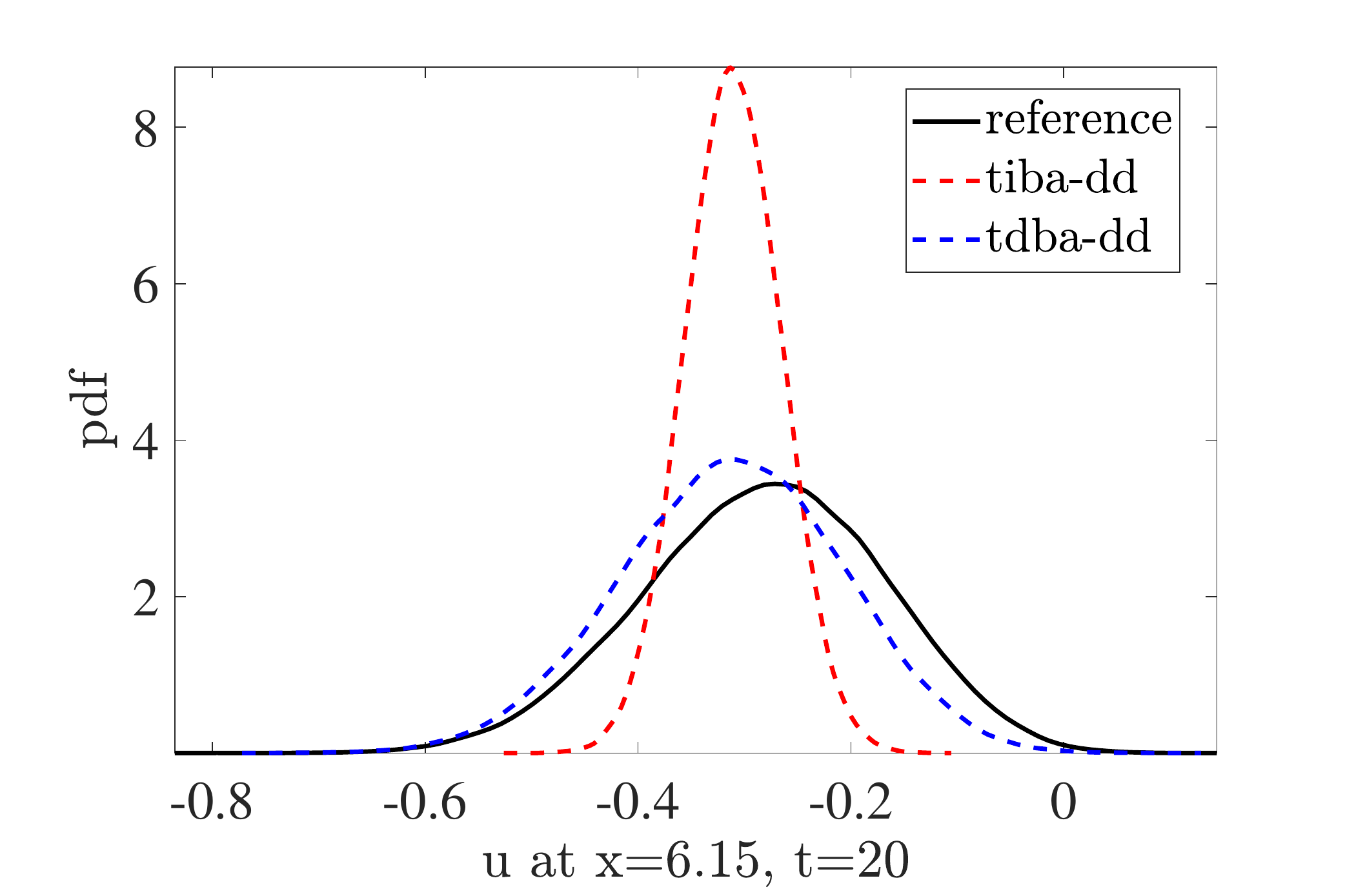}
        \caption{} \label{fig:pdf_u_at_t4001_x62_dd_xid10_p3_etad5}
    \end{subfigure}        
    \begin{subfigure}[t]{0.32\textwidth}
        \centering
        \includegraphics[scale=.25]{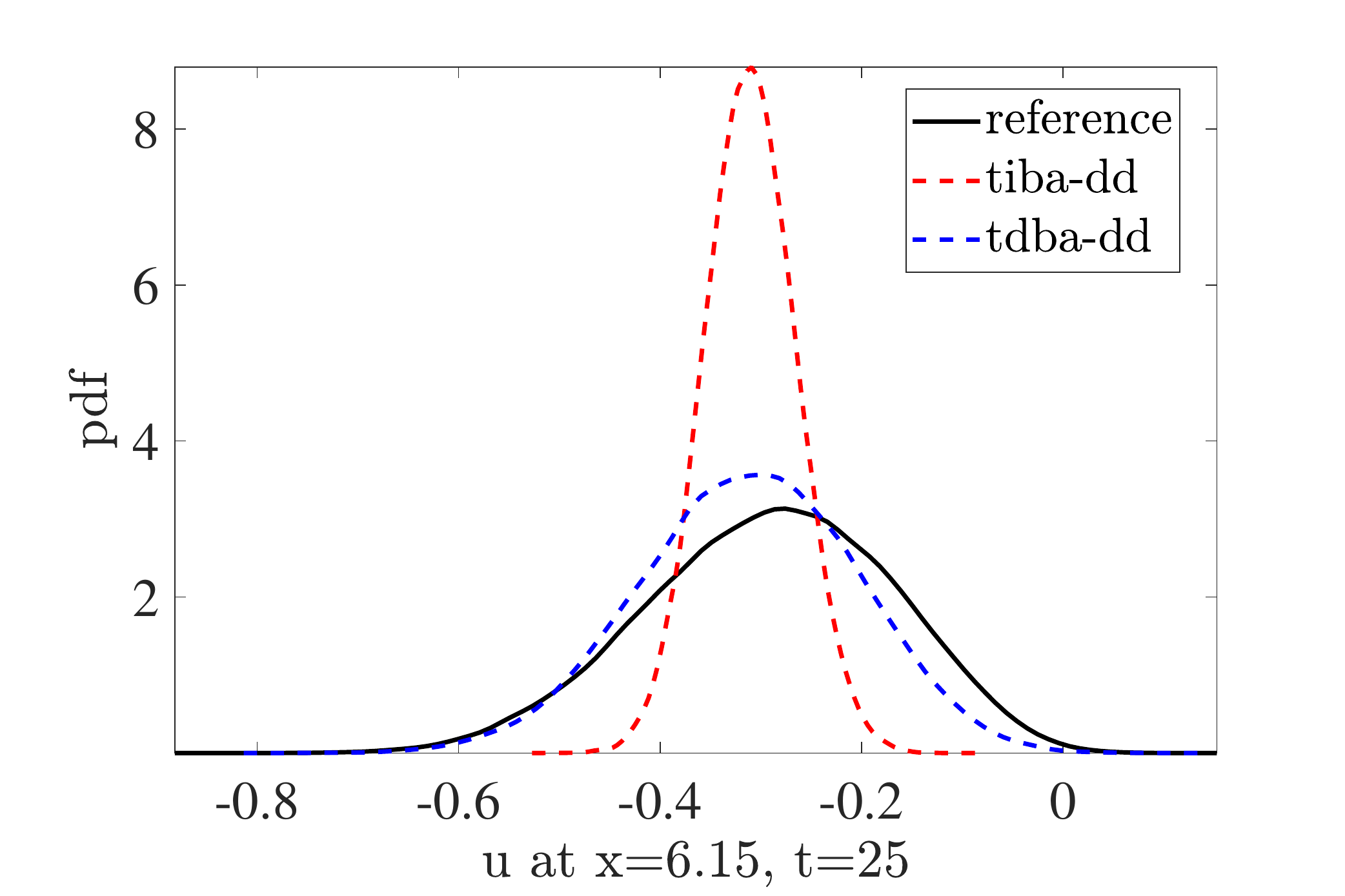}
        \caption{} \label{fig:pdf_u_at_t5001_x62_dd_xid10_p3_etad5}
    \end{subfigure} 
    \begin{subfigure}[t]{0.32\textwidth}
        \centering
        \includegraphics[scale=.25]{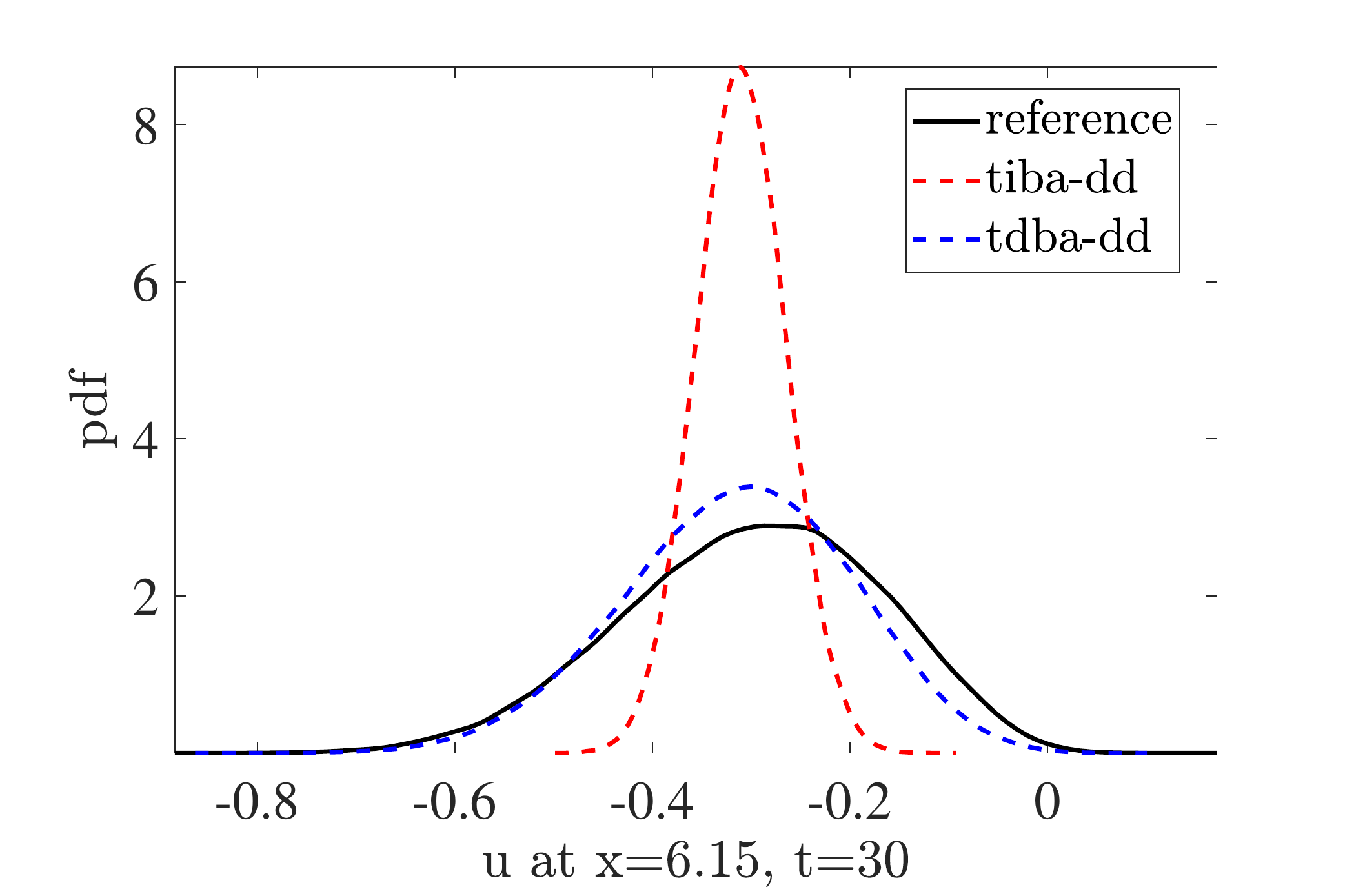}
        \caption{} \label{fig:pdf_u_at_t6001_x62_dd_xid10_p3_etad5}
    \end{subfigure}
\caption{Nonlinear diffusion equation: Probability density function of $u$ at $x=6.15$ and a). $t=0.005$, b). $t=10$, c). $t=20$, d). $t=30$, e). $t=40$ and f). $t=50.$ The reference solution (solid black) is computed with dimension $d=10$ and the reduced dimensional solutions (dashed red and blue) are computed with $r=5$.} \label{fig:pdf_u_at_t2_t6001_x62_dd_xid10_p3_etad5}
\end{figure}

Figs. \ref{fig:pdf_u_at_t2_t6001_x62_dd_xid10_p3_etad5}(a-f) show the pdf of $u$ at $x=6.15$ and $t=0.005, 10, 20, 30, 40$ and $50$ respectively.  The reference pdf (solid black line) is obtained by  solving the PDE~\eqref{eqn:1d_richards} with a $10$-dimensional random vector $\boldsymbol{\xi}$ using sparse grid level $5.$ The total number of simulations needed for this solution is $8761.$ 
In the basis adoption method with time-independent basis (the PDF is represented by the dashed red lines),  $r=5$ dimensions are used and the basis adaptation is performed only once at $t=0.005.$ The same number of dimensions ($r=5$) are used in the basis adaptation method with the time-dependent basis (PDFs are denoted by dashed blue lines). The new basis are computed at each time step prior to solving the PDE. Fig. \ref{fig:pdf_u_at_t2_t6001_x62_dd_xid10_p3_etad5} demonstrates that the pdf predicted with time-dependent basis is significantly more accurate than the PDF obtained with fixed basis.

\begin{figure}[ht!]
    \centering
    \begin{subfigure}[t]{0.48\textwidth}
        \centering
        \includegraphics[scale=.23]{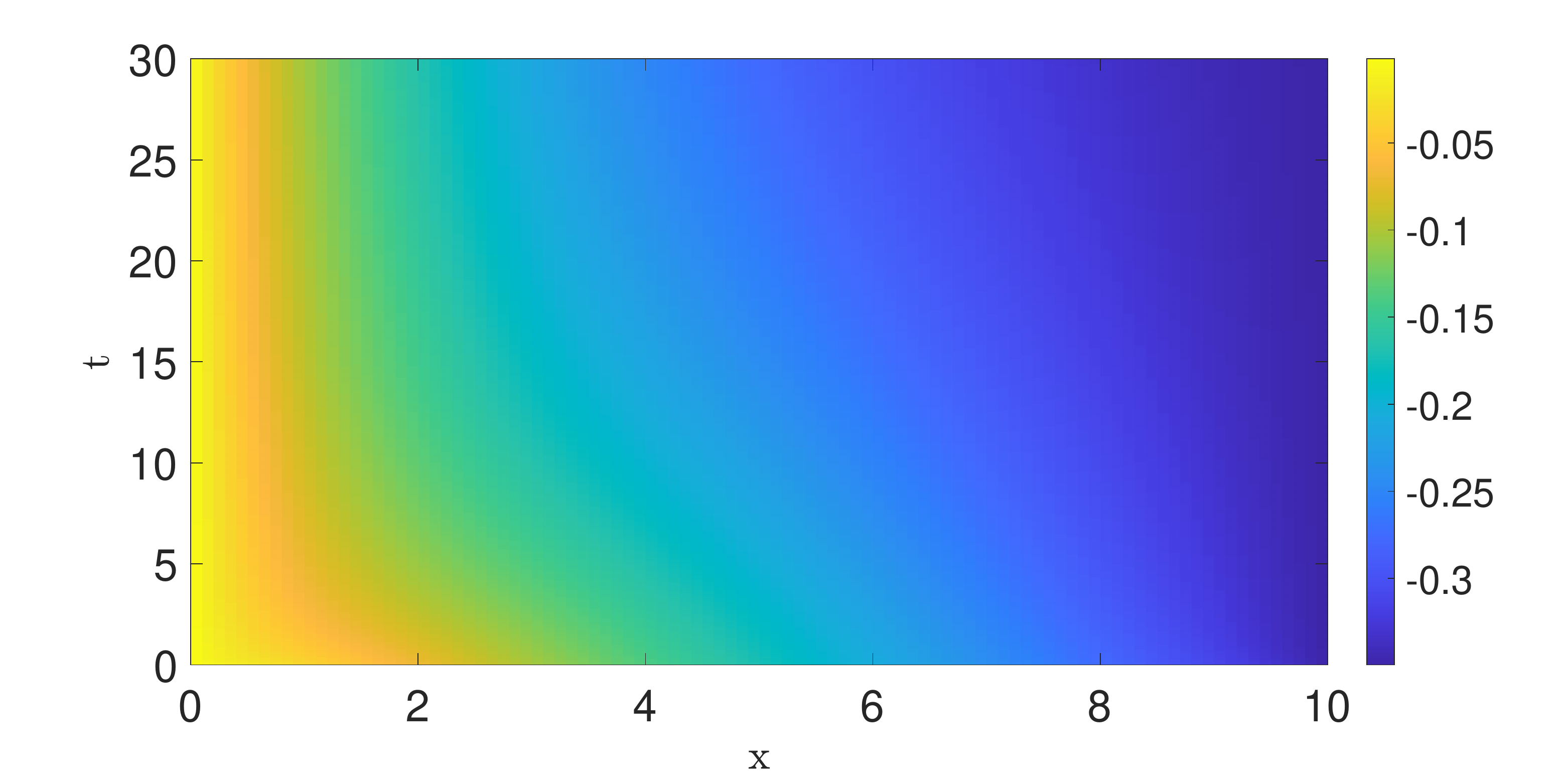}
        \caption{Mean} \label{fig:Yt_mean_xid10}
    \end{subfigure}        
    \begin{subfigure}[t]{0.48\textwidth}
        \centering
        \includegraphics[scale=.23]{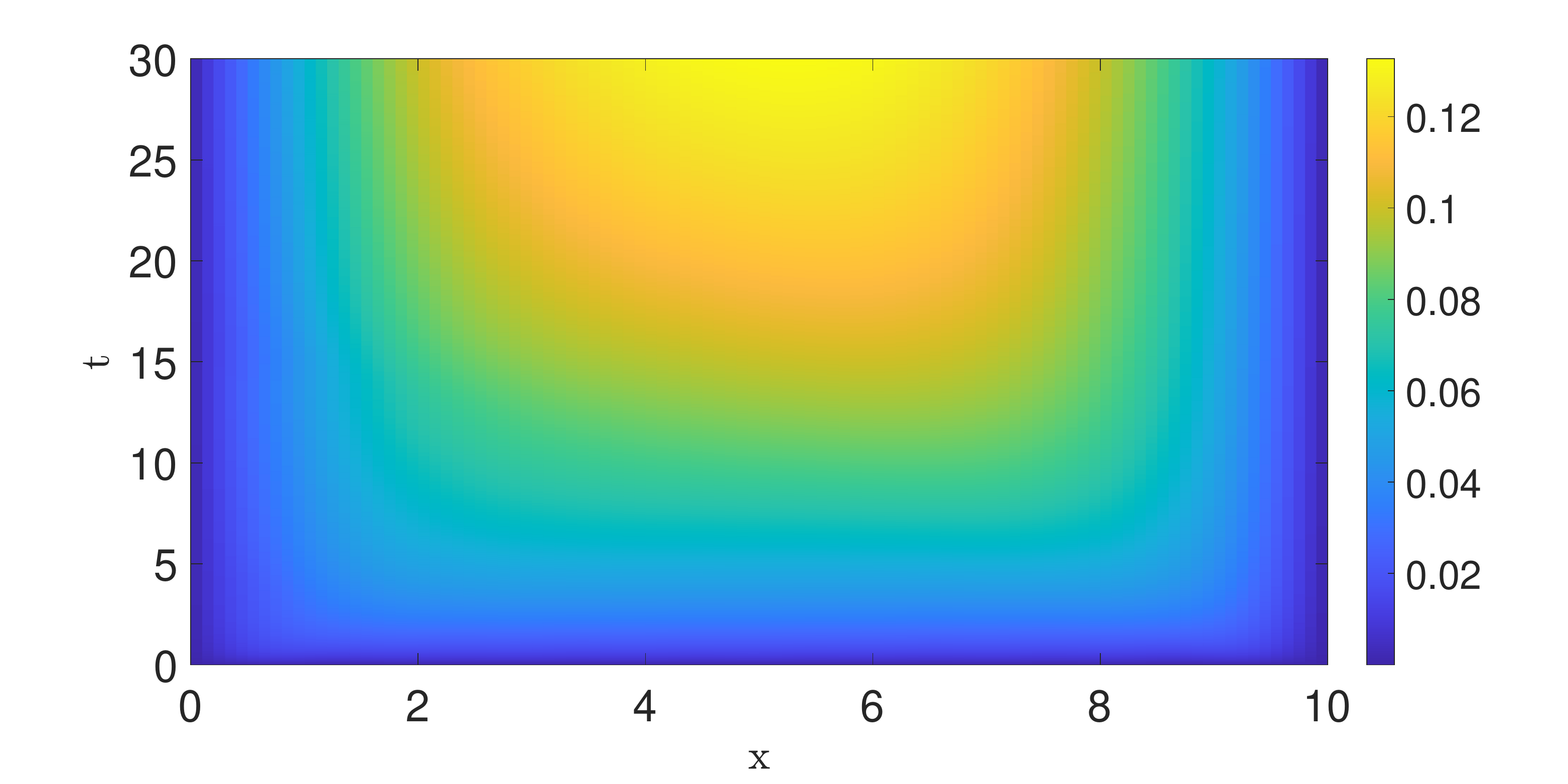}
        \caption{Standard deviation} \label{fig:Yt_sdev_xid10}
    \end{subfigure} 
    \begin{subfigure}[t]{0.48\textwidth}
        \centering
        \includegraphics[scale=.23]{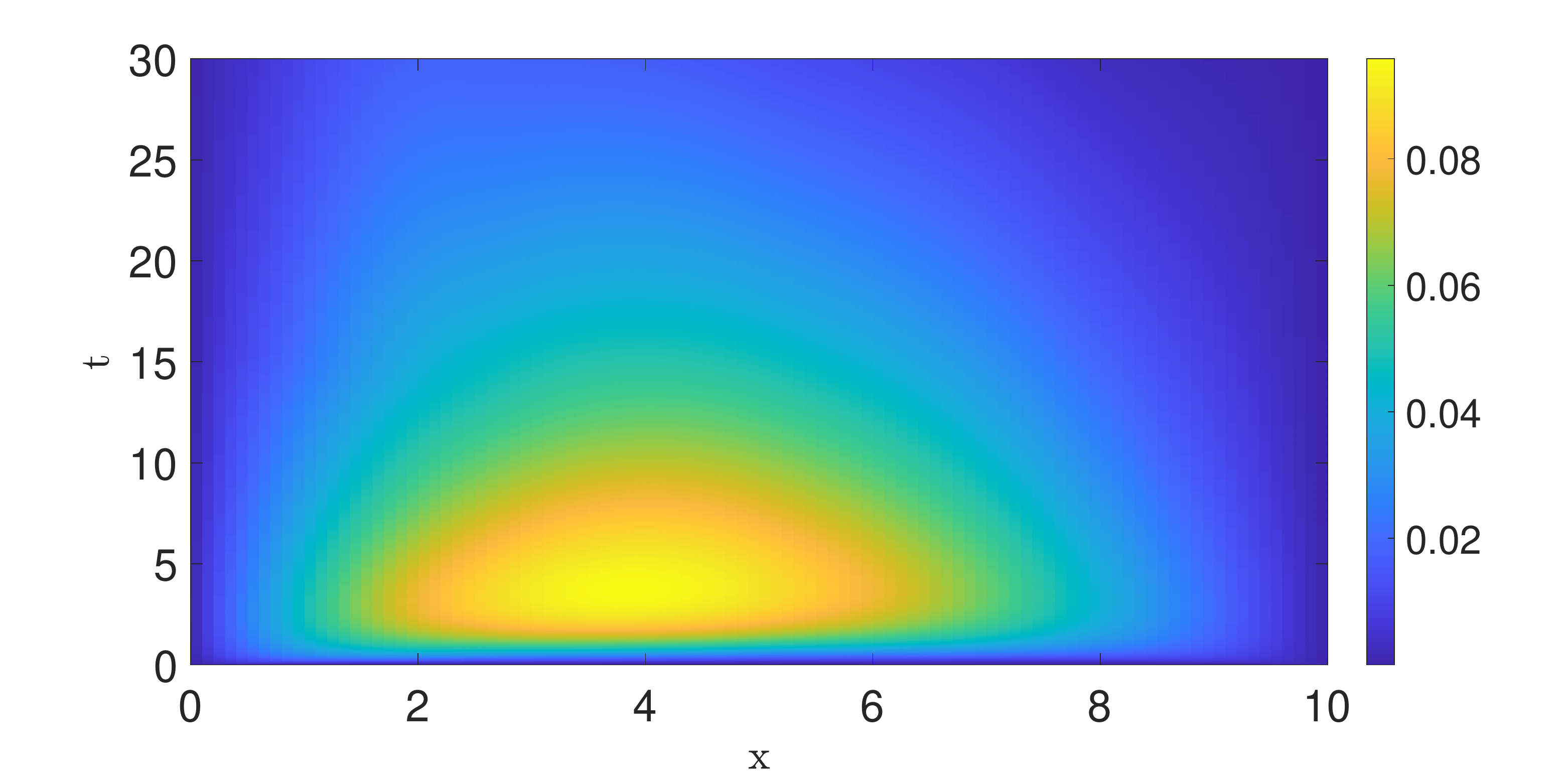}
        \caption{Abs. point error in mean} \label{fig:Yt_DD_mean_etad5_ti_error}
    \end{subfigure}        
    \begin{subfigure}[t]{0.48\textwidth}
        \centering
        \includegraphics[scale=.23]{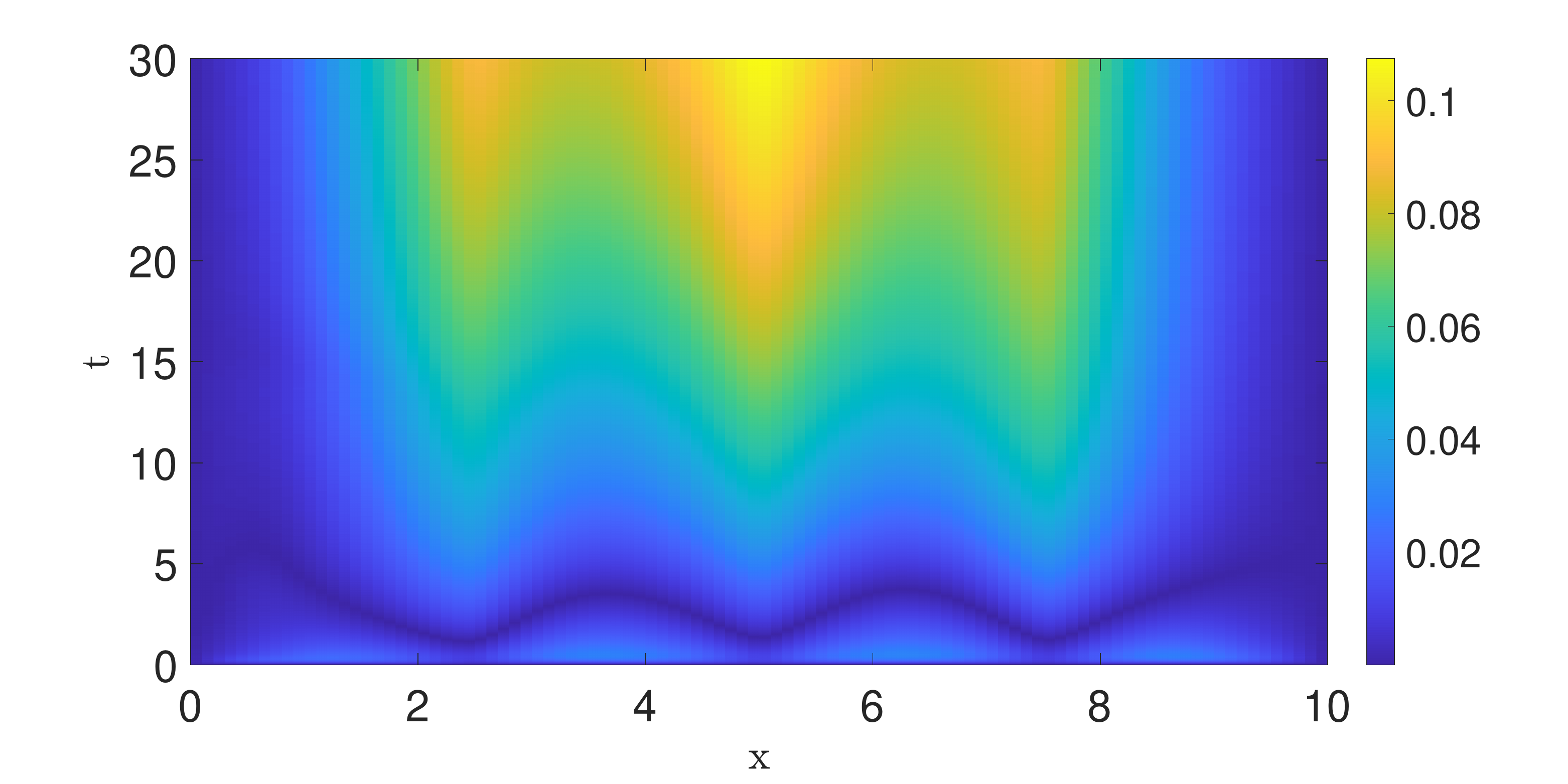}
        \caption{Abs. point error in standard deviation} \label{fig:Yt_DD_sdev_etad5_ti_error}
    \end{subfigure} 
    \begin{subfigure}[t]{0.48\textwidth}
        \centering
        \includegraphics[scale=.23]{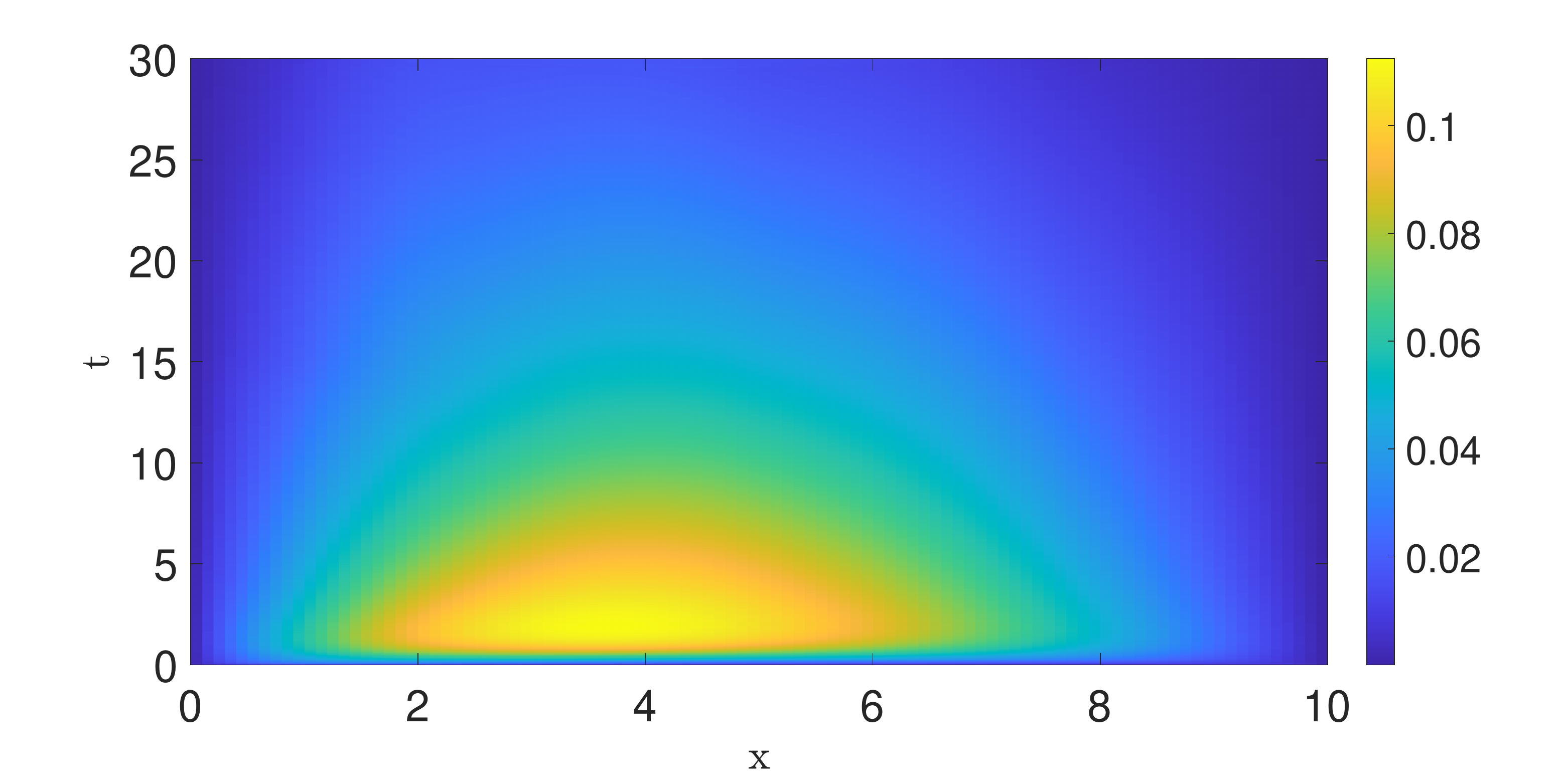}
        \caption{Abs. point error in mean} \label{fig:Yt_DD_mean_etad5_tdba_error}
    \end{subfigure}        
    \begin{subfigure}[t]{0.48\textwidth}
        \centering
        \includegraphics[scale=.23]{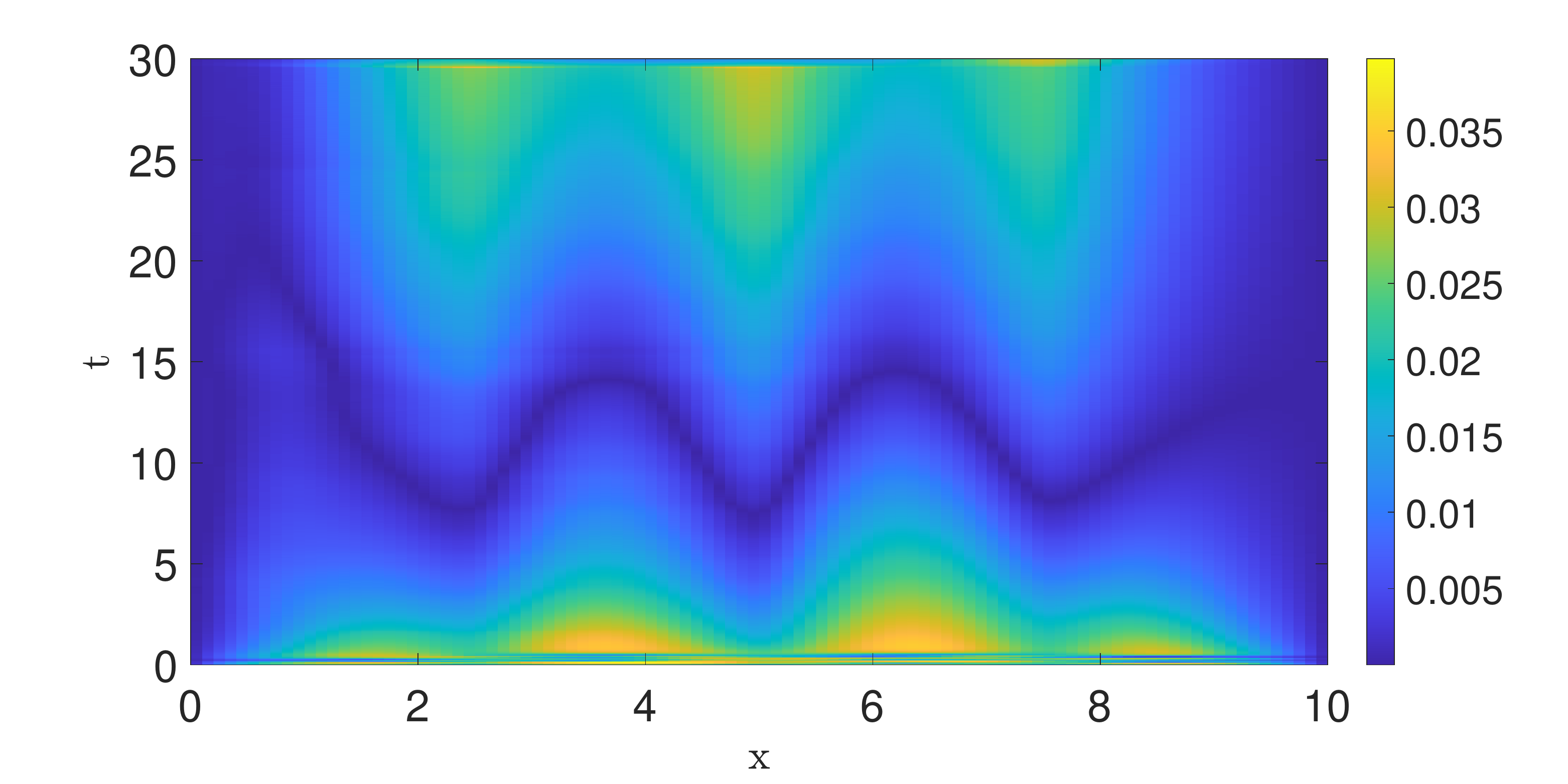}
        \caption{Abs. point error in standard deviation} \label{fig:Yt_DD_sdev_etad5_tdba_error}
    \end{subfigure} 
\caption{Richards equation: a) mean of reference solution ($d=10$) and b) standard deviation of reference solution ($d=10$) c) error in mean of low dimensional solution ($r=5$) with fixed basis adaptation and d) error in standard deviation of low dimensional solution ($r=5$) with fixed basis adaptation e) error in mean of low dimensional solution ($r=5$) with time-dependent basis adaptation and f) error in standard deviation of low dimensional solution ($r=5$) with time-dependent basis adaptation as a function of $x$ (horizontal axis) and $t$ (vertical axis).} \label{fig:Yt_DD_sdev_etad5_error}
\end{figure}

Figs. \ref{fig:Yt_mean_xid10} and \ref{fig:Yt_sdev_xid10} show the the mean and standard deviation of the reference solution as a function of $x$ (horizontal axis) and $t$ (vertical axis). Figs. \ref{fig:Yt_DD_mean_etad5_ti_error} and \ref{fig:Yt_DD_sdev_etad5_ti_error} present the absolute point error in the mean and standard deviation of the solution using fixed basis adaptation. Figs. \ref{fig:Yt_DD_mean_etad5_tdba_error} and \ref{fig:Yt_DD_sdev_etad5_tdba_error} display the absolute point error in the mean and standard deviation of the solution of the Richards using time-dependent basis adaptation compared to the reference solution. While the point errors of the mean are similar in the  the  time-independent and time-dependent basis adaptation methods, the point errors in the standard deviation are approximately 50\% smaller with the time-dependent basis. 

In the above examples, the solutions are obtained with the domain decomposition method by dividing the domain into $4$ non-overlapping subdomains and in each subdomain the solution is computed using $781$ simulations corresponding to the reduced dimension $r=5$ and sparse-grid level $5.$ The total number of simulations is $781 \times 4 = 3124$, which is significantly smaller than $8761$ simulations in the reference solution.  The iterative method for computing values at the interface between the subdomains does not add significantly to the computational cost.

\begin{figure}[ht!]
    \centering
    \begin{subfigure}[t]{0.48\textwidth}
        \centering
        \includegraphics[scale=.23]{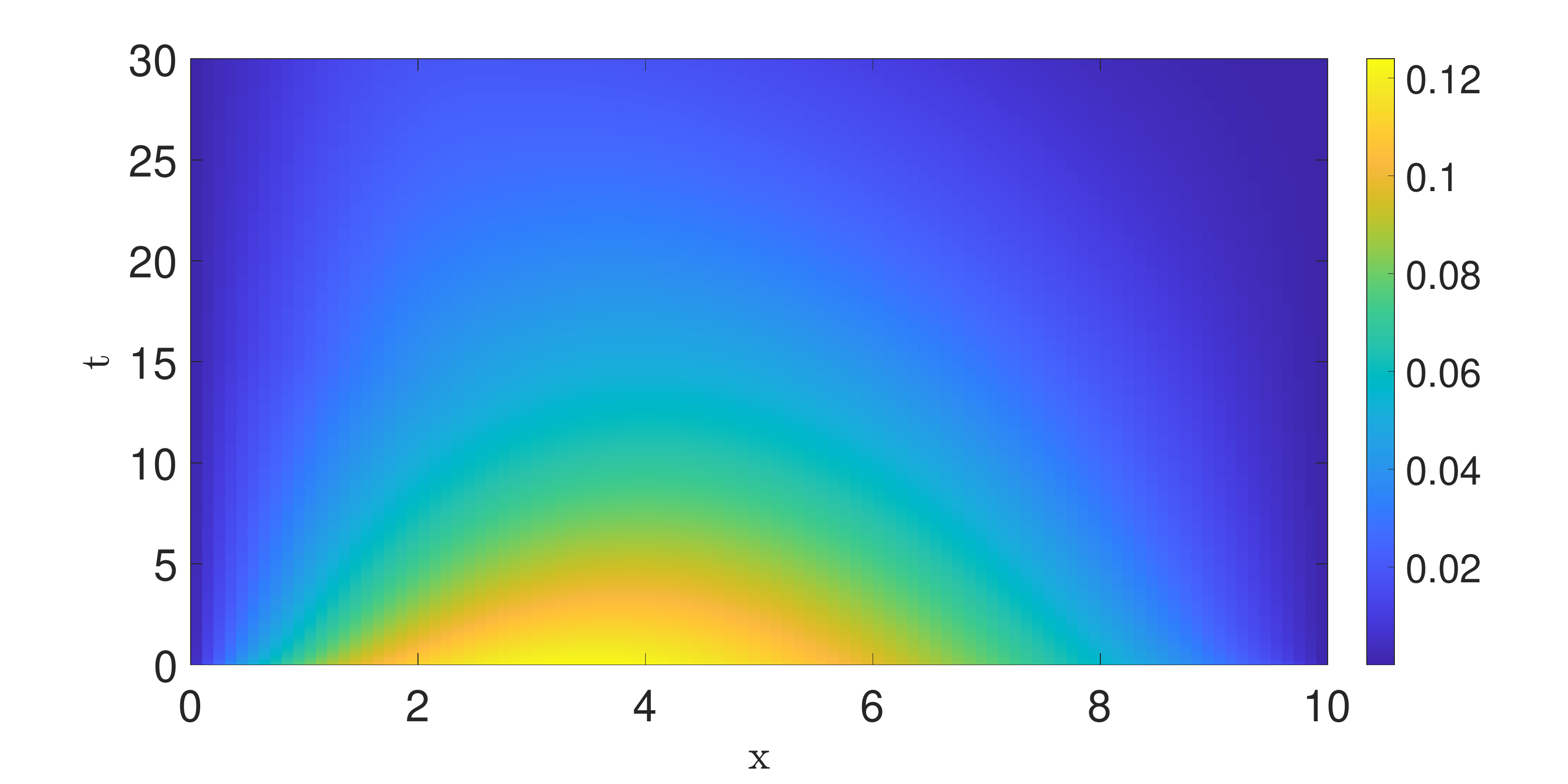}
        \caption{Abs. point error in mean} \label{fig:Yt_wo_DD_mean_etad5_ti_error}
    \end{subfigure}        
    \begin{subfigure}[t]{0.48\textwidth}
        \centering
        \includegraphics[scale=.23]{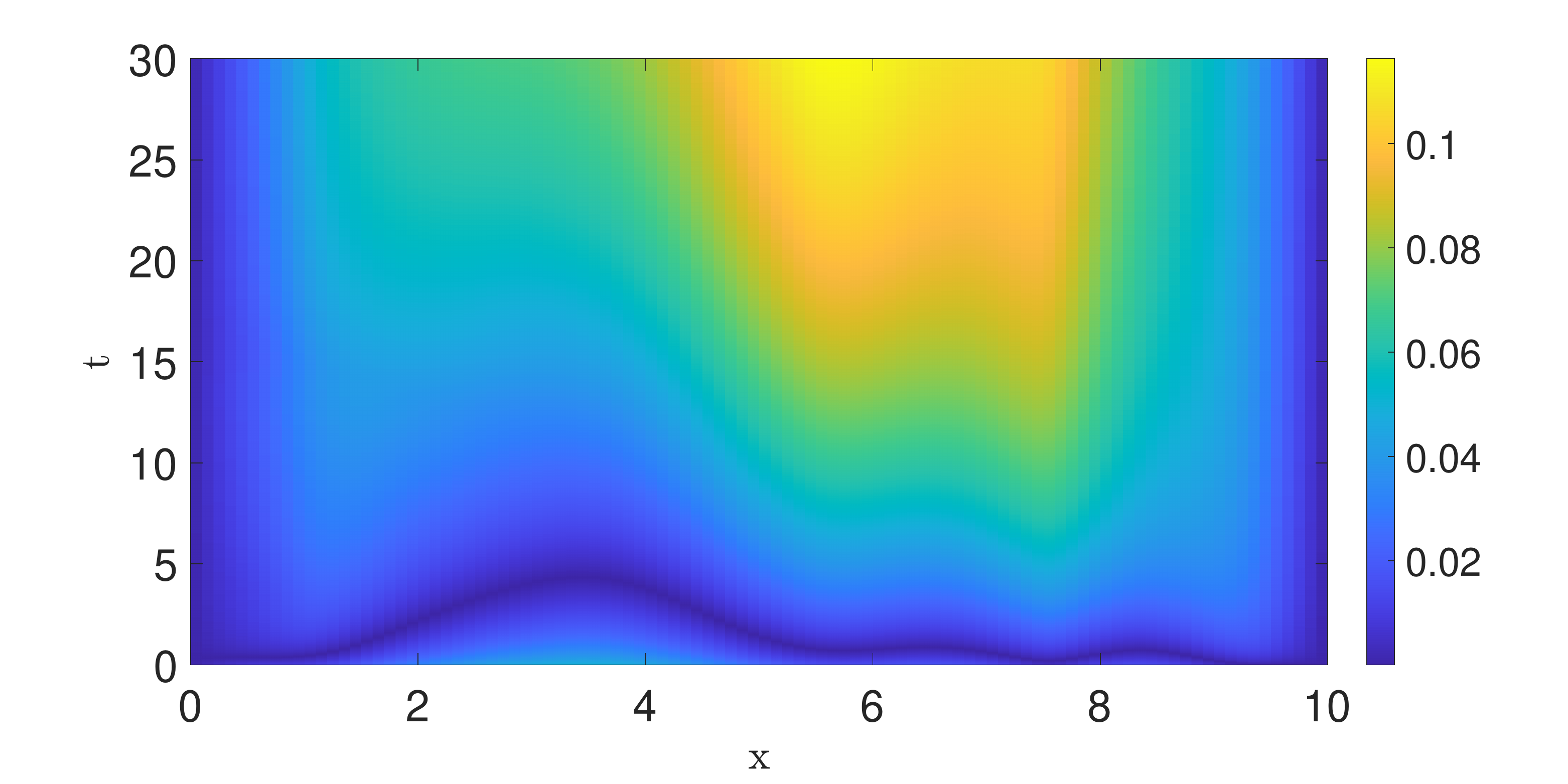}
        \caption{Abs. point error in standard deviation} \label{fig:Yt_wo_DD_sdev_etad5_ti_error}
    \end{subfigure} 
    \begin{subfigure}[t]{0.48\textwidth}
        \centering
        \includegraphics[scale=.23]{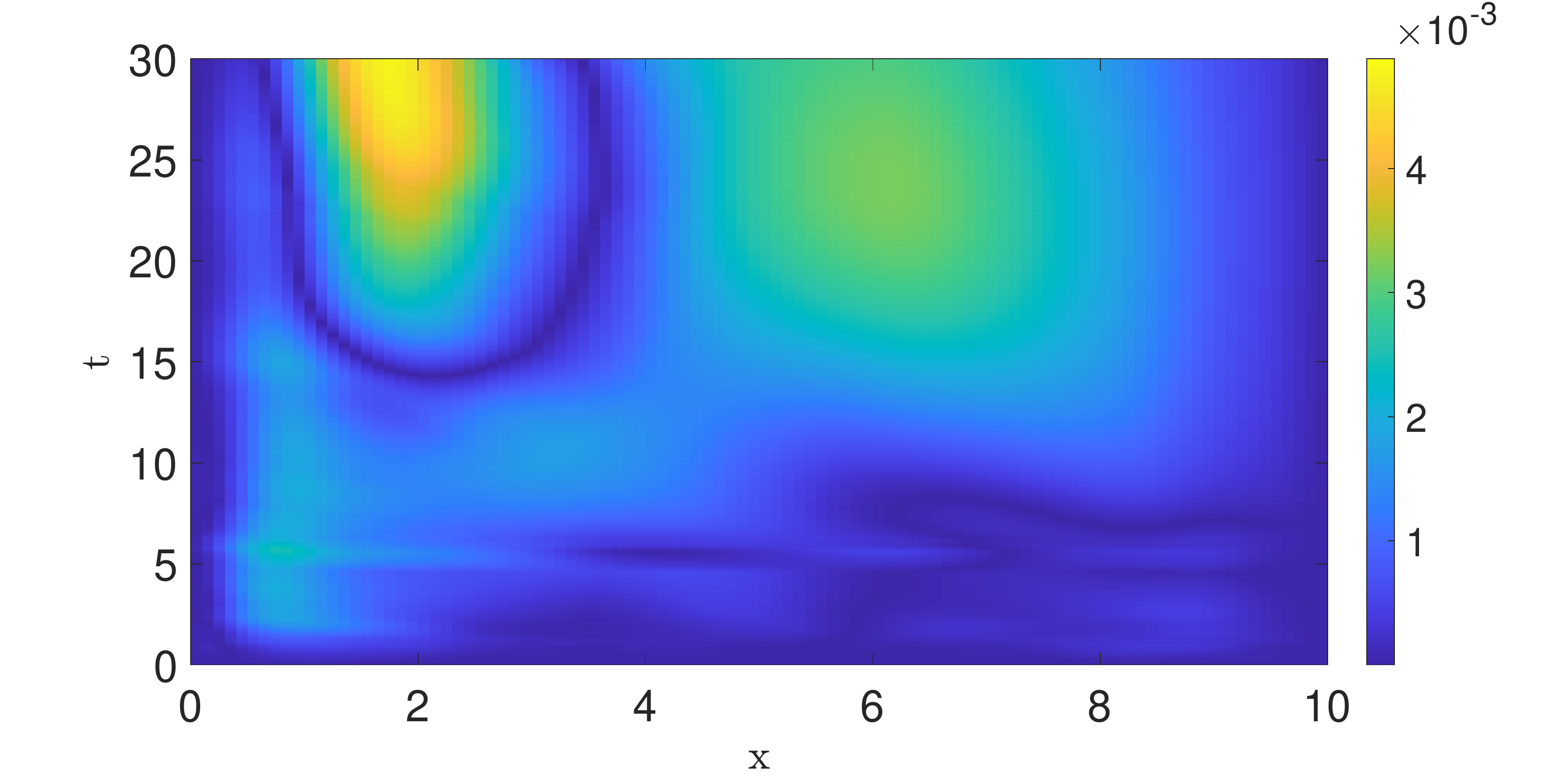}
        \caption{Abs. point error in mean} \label{fig:Yt_wo_DD_mean_etad5_tdba_error}
    \end{subfigure}        
    \begin{subfigure}[t]{0.48\textwidth}
        \centering
        \includegraphics[scale=.23]{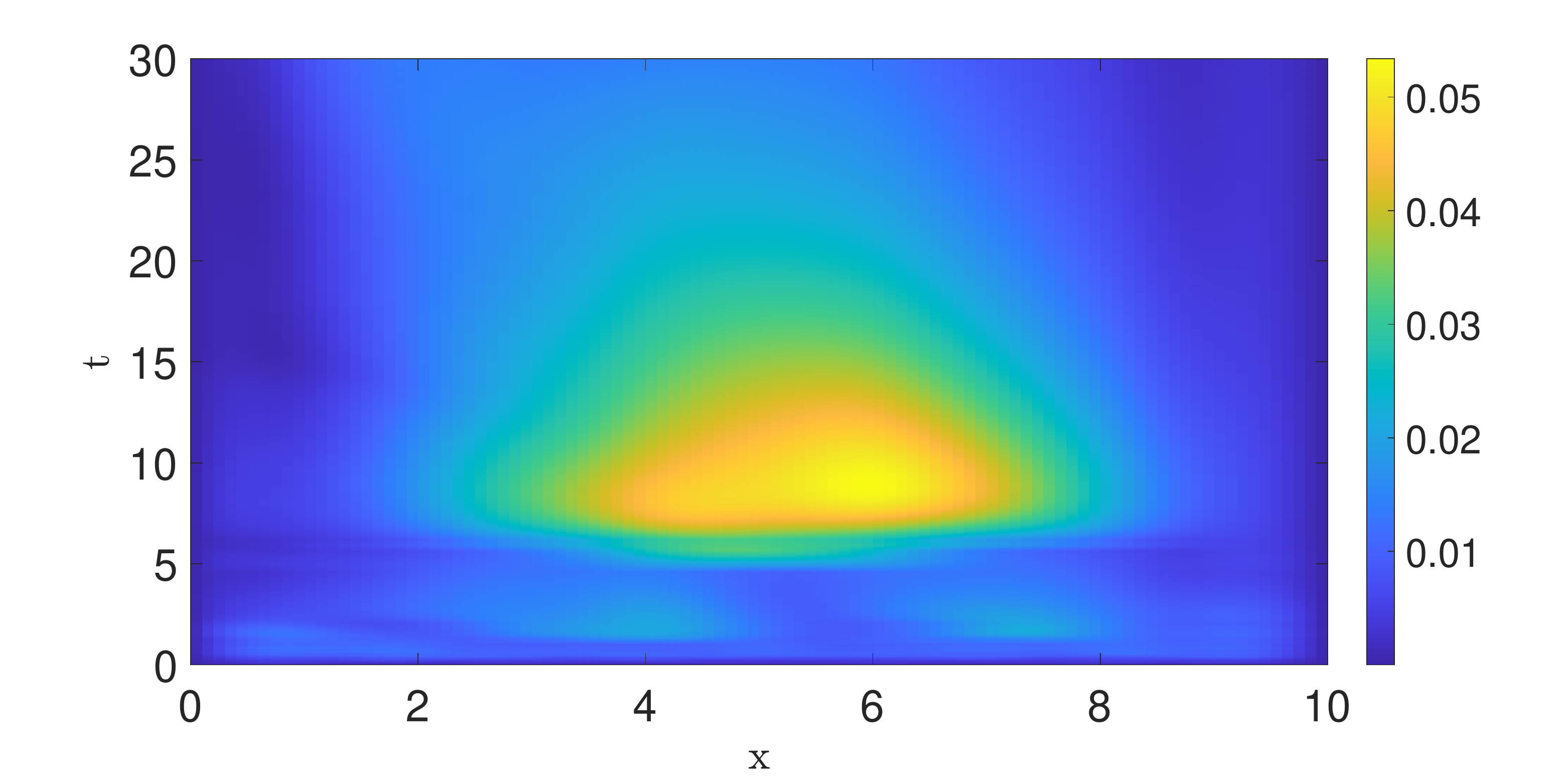}
        \caption{Abs. point error in standard deviation} \label{fig:Yt_wo_DD_sdev_etad5_tdba_error}
    \end{subfigure} 
\caption{Richards equation: Basis adaptation without domain decomposition a) error in mean of low dimensional solution ($r=5$) with fixed basis adaptation, b) error in standard deviation of low dimensional solution ($r=5$) with fixed basis adaptation c) error in mean of low dimensional solution ($r=5$) with time-dependent basis adaptation and d) error in standard deviation of low dimensional solution ($r=5$) with time-dependent basis adaptation as a function of $x$ (horizontal axis) and $t$ (vertical axis).} \label{fig:Yt_wo_DD_sdev_etad5_error}
\end{figure}

Finally, in Fig. \ref{fig:Yt_wo_DD_sdev_etad5_error} we present errors in the solution obtained with the time-dependent basis adaptation without domain decomposition. These errors are smaller than in the time-independent basis adaptation method but larger than in the time-dependent basis adaptation method with the domain decomposition.

\section{Conclusions} \label{sec:conclusions}
We proposed a dimension reduction method with time-dependent basis adaption and the spatial domain decomposition for stochastic PDEs. The time-dependent basis adaptation is computed using an inexpensive calculation of the full dimensional stochastic solution in the whole spatial domain using lower sparse-grid level (2 in this work).

We have used the proposed approach for solving the time-dependent one-dimensional linear and nonlinear diffusion equations with random diffusion coefficient. 
Our results show that the proposed approach significantly outperforms  the time-independent basis adaptation method. Although the time-dependent basis adaptation method alone is quite effective, combining with the domain decomposition method further improved the results. To the best of our knowledge, this is the first time-dependent dimension reduction approach for time-dependent stochastic PDEs. In our future work, we will extend our approach to two and three dimensional stochastic PDEs.     

\section{Acknowledgments}
This research was supported by the U.S. Department of Energy, Office of Science, Office of Advanced 
Scientific Computing Research as part of the ``Uncertainty Quantification For Complex Systems Described by Stochastic Partial Differential Equations'' project.  Pacific Northwest National Laboratory is operated by Battelle for the DOE under Contract DE-AC05-76RL01830.

\bibliographystyle{elsarticle-num} 

\bibliography{references}

\begin{thebibliography}{10}
\expandafter\ifx\csname url\endcsname\relax
  \def\url#1{\texttt{#1}}\fi
\expandafter\ifx\csname urlprefix\endcsname\relax\def\urlprefix{URL }\fi
\expandafter\ifx\csname href\endcsname\relax
  \def\href#1#2{#2} \def\path#1{#1}\fi

\bibitem{tipireddy2014basis}
R.~Tipireddy, R.~Ghanem, Basis adaptation in homogeneous chaos spaces, Journal
  of Computational Physics 259 (2014) 304--317.

\bibitem{Tipireddy2013}
R.~Tipireddy, {Algorithms for stochastic Galerkin projections: Solvers, basis
  adaptation and multiscale modeling and reduction}, {T}heses, {University of
  Southern California} (Aug. 2013).

\bibitem{Tsilifis2016}
P.~Tsilifis, R.~Ghanem, Reduced {Wiener Chaos} representation of random fields
  via basis adaptation and projection, arXiv:1603.04803v3.

\bibitem{Constantine2014}
P.~G. Constantine, E.~Dow, Q.~Wang, Active subspace methods in theory and
  practice: Applications to kriging surfaces, Methods and Algorithms for
  Scientific Computing 36 (2014) A1500--A1524.

\bibitem{Li2016}
W.~Li, G.~Lin, B.~Li, Inverse regression-based uncertainty quantification
  algorithms for high-dimensional models: Theory and practice, Journal of
  Computational Physics 321 (2016) 259--278.

\bibitem{tipireddy2017basis}
R.~Tipireddy, P.~Stinis, A.~M. Tartakovsky, Basis adaptation and domain
  decomposition for steady-state partial differential equations with random
  coefficients, Journal of Computational Physics 351 (2017) 203--215.

\bibitem{tipireddy2018stochastic}
R.~Tipireddy, P.~Stinis, A.~Tartakovsky, Stochastic basis adaptation and
  spatial domain decomposition for partial differential equations with random
  coefficients, SIAM/ASA Journal on Uncertainty Quantification 6~(1) (2018)
  273--301.

\bibitem{Chen2015}
Y.~Chen, J.~Jakeman, C.~Gittelson, D.~Xiu, Local polynomial chaos expansion for
  linear differential equations with high dimensional random inputs, {SIAM}
  Journal on Scientific Computing 37~(1) (2015) A79--A102.

\bibitem{Toselli2005}
A.~Toselli, O.~B. Widlund, Domain Decomposition Methods--Algorithms and Theory,
  Springer-Verlag, 2005.

\bibitem{Doostan2007}
A.~Doostan, R.~G. Ghanem, J.~Red-Horse, Stochastic model reduction for chaos
  representations, Computer Methods in Applied Mechanics and Engineering 196
  (2007) 3951--3966.

\bibitem{Levy1999}
A.~Levy, J.~Rubinstein, {H}ilbert--space {K}arhunen--{L}o\'eve transform with
  application to image analysis, Journal of The Optical Society of America A 16
  (1999) 28--35.

\bibitem{Kirby1992}
M.~Kirby, Minimal dynamical systems from {PDE}s using {S}obolev eigenfunctions,
  Physica D: Nonlinear Phenomena 57 (1992) 466--475.

\bibitem{Silverman1996}
B.~W. Silverman, Smoothed functional principal components analysis by choice of
  norm, The Annals of Statistics 24 (1996) 1--24.

\bibitem{Berkooz1993539}
G.~Berkooz, P.~Holmes, J.~Lumley, The proper orthogonal decomposition in the
  analysis of turbulent flows, Annual Review of Fluid Mechanics 25~(1) (1993)
  539--575.

\bibitem{Christensen1999}
E.~A. Christensen., M.~Br{\o}ns, J.~N. S{\o}rensen, Evaluation of proper
  orthogonal decomposition--based decomposition techniques applied to
  parameter-dependent nonturbulent flows, {SIAM} Journal on Scientific
  Computing 21~(4) (1999) 1419--1434.

\bibitem{pan1995transformed}
L.~Pan, P.~J. Wierenga, A transformed pressure head-based approach to solve
  richards' equation for variably saturated soils, Water Resources Research
  31~(4) (1995) 925--931.

\end{thebibliography}

\end{document}